\documentclass{article}
 \usepackage{pict2e}
\RequirePackage{amssymb,amsthm,amsmath}

\newtheorem{Theorem}{Theorem}[section]
\newtheorem{Corollary}[Theorem]{Corollary}
\newtheorem{Proposition}[Theorem]{Proposition}
\newtheorem{Lemma}[Theorem]{Lemma}
\newtheorem{Conjecture}[Theorem]{Conjecture}

 \renewcommand{\stretch}{\mathrm{stretch}}
 \newcommand{\len}{\, \mathrm{len}}
 \newcommand{\ST}{\mathrm{ST}}
\newcommand{\cworst}{c_{\mbox{{\tiny worst}}}}
\newcommand{\cave}{c_{\mbox{{\tiny ave}}}}
\def\ind{{\rm 1\hspace{-0.90ex}1}}
 \newcommand{\eps}{\varepsilon}
\newcommand{\Psave}{\Psi^{\mbox{\tiny ave}}}
\newcommand{\Psworst}{\Psi^{\mbox{\tiny worst}}}
\newcommand{\bz}{{\mathbf z}}
\newcommand{\bt}{{\mathbf t}}
\newcommand{\bZ}{{\mathbf Z}}
\newcommand{\NN}{\mbox{${\mathcal N}$}}
\newcommand{\sfrac}[2]{{\textstyle\frac{#1}{#2}}}
 \newcommand{\Ex}{{\mathbb{E}}}
 \renewcommand{\Pr}{{\mathbb{P}}}
 \newcommand{\ave}{\mathrm{ave}}
 \newcommand{\Reals}{{\mathbb{R}}}
\newcommand{\pmut}{p_{\mathrm{mut}}}
\newcommand{\pnot}{p_{\mathrm{not}}}
\newcommand{\origin}{\mathbf{0}}
\newcommand{\area}{\mathrm{area}}
 \newcommand{\intersect}{\mathrm{intersect}}
 \newcommand{\cone}{\mathrm{cone}}
\begin{document}

\title{The Stretch - Length Tradeoff in Geometric Networks: Average Case and Worst Case Study}

\author{David Aldous\\Department of Statistics\\U.C. Berkeley CA 94720 \and Tamar Lando\\Department of Philosophy\\Columbia University\\1150 Amsterdam Ave\\New York NY 10027}

\maketitle

\begin{abstract}
Consider a network linking the points of a rate-$1$ Poisson point process on the plane.
Write $\Psave(s)$ for the minimum possible mean length per unit area of such a network,  subject to the constraint that the route-length between every pair of points is at most $s$ times the Euclidean distance.  We give upper and lower bounds on the function  $\Psave(s)$, and on the analogous 
``worst-case" function $\Psworst(s)$ where the point configuration is arbitrary subject to average density one per unit area.  
Our bounds are numerically crude, but raise the question of whether 
there is an exponent $\alpha$ such that each function has 
$\Psi(s) \asymp (s-1)^{-\alpha}$ as $s \downarrow 1$.
\end{abstract}

\section{Introduction}\label{sec:1}
The topic 
 {\em geometric spanner networks} \cite{MR2289615}
 concerns design of networks on arbitrary sets of vertices in the plane (or higher dimensions).
The interpretation of ``size" of the network is sometimes as number of edges
and sometimes as network length (sum of Euclidean edge lengths).  
Similarly, 
the interpretation of within-network distance between two vertices $v,w$ is sometimes taken as minimum number of edges 
of a route between them ({\em hop length}) 
and sometimes as shortest total length (sum of Euclidean edge lengths) of a route between them 
({\em route length} $r(v,w)$).   
In the latter setting, 
how well the network provides short routes 
is often measured by a statistic such as
\begin{equation}
S := 
\max_{v \neq w} \frac{r(v,w)}
{d(v,w)}
\ge 1  
\label{S-def}
\end{equation}
 where $d(v,w)$ denotes straight line (Euclidean) distance. 
The statistic $S$ is called the {\em stretch} 
or {\em spanning ratio} of the network, and 
a network with stretch $S$ is called an $S$-{\em spanner}.

Most work on this topic has emphasized algorithms -- either algorithms for constructing spanners, 
or the use of spanners in algorithms for computational geometry problems.
We address a more fundamental geometric question: what is the tradeoff between stretch and network length? 
In formulating a mathematical question we have in mind 
the example of an inter-city road network (rather than, say, a wireless communication network) and indeed we find it helpful to use the vivid natural language of {\em cities, roads, junctions} in place of the mathematical language of {\em vertices, edges, Steiner points}.  

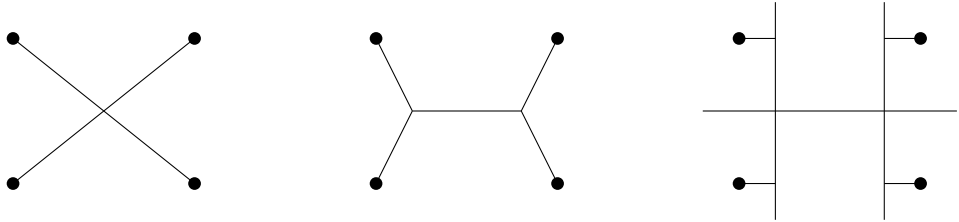
\begin{figure}
 \setlength{\unitlength}{0.095in}
 \begin{picture}(52,15)(0,-3)
\put(0,0){\circle*{0.7}}
\put(10,0){\circle*{0.7}}
\put(0,8){\circle*{0.7}}
\put(10,8){\circle*{0.7}}
\put(0,0){\line(5,4){10}}
\put(10,0){\line(-5,4){10}}
\put(20,0){\circle*{0.7}}
\put(30,0){\circle*{0.7}}
\put(20,8){\circle*{0.7}}
\put(30,8){\circle*{0.7}}
\put(20,0){\line(1,2){2}}
\put(20,8){\line(1,-2){2}}
\put(30,0){\line(-1,2){2}}
\put(30,8){\line(-1,-2){2}}
\put(22,4){\line(1,0){6}}
\put(40,0){\circle*{0.7}}
\put(50,0){\circle*{0.7}}
\put(40,8){\circle*{0.7}}
\put(50,8){\circle*{0.7}}
\put(42,-2){\line(0,1){12}}
\put(48,-2){\line(0,1){12}}
\put(38,4){\line(1,0){14}}
\put(40,0){\line(1,0){2}}
\put(40,8){\line(1,0){2}}
\put(50,0){\line(-1,0){2}}
\put(50,8){\line(-1,0){2}}
\end{picture}
\caption{Illustration of possible networks. 
The left diagram shows a network on $4$ cities which is connected 
(a driver can switch roads at the junction where they cross), the center diagram shows that other junctions  may be created, and the right diagram 
(envisaged as part of a larger network) shows that roads need not be closely related to cities at all.}
\end{figure}

Our notion of ``network" is the general notion suggested by real-world  road networks, illustrated in Fig. 1.  The only implicit convention is that networks be connected and consist of line segments (rather than curves).  
Note this convention differs from the more familiar and more restrictive
 assumption that an edge can only be the line segment $(v,w)$ 
between two of the given vertices.  
There does not seem to be standard terminology to emphasize the distinction: we will write {\em Steiner network} for our setting and 
{\em graph network} for the more restrictive
 assumption.

Our underlying setting is
a configuration of $n$ cities at arbitrary positions 
$\bz_n = (z_1, \ldots, z_n)$ 
in a square of area $n$.  
For a network $\NN$ connecting these cities, write 
$S(\NN)$ for the statistic (\ref{S-def}) and write 
\[ L(\NN) = \sfrac{1}{n} \times \mbox{(network length of $\NN$)} \]
for normalized network length.
We then define

\[
\psi_n(\bz_n,s) := \inf \{L(\NN): \ S(\NN) \leq  s \}  \] 
the infimum over all networks $\NN$ connecting the cities $\bz_n$.
So this quantifies the optimal trade-off between length and stretch for the given configuration.  
We can now consider in parallel the worst-case, that is 
$\sup_{\bz_n} \psi_n(\bz_n,s) $, 
and the average case 
$\Ex \psi_n(\bZ_n,s) $, 
where $\bZ_n$ consists of $n$ independent uniform random positions in the area-$n$ square.   
The purpose of this set-up is that it is intuitively obvious  that there must exist limit functions 
\[ \Psworst(s) = \lim_{n \to \infty} \sup_{\bz_n} \psi_n(\bz_n,s)  \]
\[ \Psave(s) =  \lim_{n \to \infty}  \Ex \psi_n(\bZ_n,s) \]
where 
$0 < \Psave(s) \leq \Psworst(s) \leq \infty$ for $1<s<\infty$. 

The goal of this paper is to study the functions $\Psave$ and $\Psworst$. 
As probabilists the authors are primarily interested in the average-case setting, but it seems natural to treat the worst-case setting in parallel.
Here is what we shall do.

\begin{itemize}
\item Prove existence of these limit functions (section \ref{sec-2}).
\item Prove that their $s \to \infty$ limits are equal to (rather than greater than) the associated Steiner tree constants (section \ref{sec-3}).
\item See what upper bounds can be derived from known results 
for graph networks (sections \ref{sec:already} and \ref{sec-5}).
\item Derive upper bounds on $\Psworst(s)$ from elementary constructions where one first lays down a regular 
network of roads without paying attention to city positions,
and then adds local links from 
cities to the network (section \ref{sec-4}).
\item Derive upper bounds on $\Psave(s)$  from Steiner network analogs 
of the $\theta$-graphs discussed in section \ref{sec-5} (section \ref{sec-99}).
\item Derive lower bounds on $\Psave(s)$ for small $s$, based on the stochastic geometry relationship between network length and rate of intersections with a typical line
(section \ref{sec-6}). 
\item Derive lower bounds on $\Psworst(s)$ based on a notion of ``local optimality" for specific networks on specific configurations (section \ref{sec-7}). 
\end{itemize}
The sections are, to a large extent, independent of each other. 
Our description of $\Psave(s)$ as an $n \to \infty$ limit was intended to 
facilitate comparsion with the worst-case setting.
A more abstract interpretation of $\Psave(s)$ in terms of networks on 
a Poisson point process on $\Reals^2$ is given in section \ref{sec-Pois-interpret} .

Getting explicit values for these functions analytically  seems impossible.  Even getting convicing numerical values would provide a 
 challenge for designers of heuristic algorithms, and we have not attempted to do so.   
Our bounds are numerically crude; this paper is intended to 
initiate study of these functions, not to give definitive results.

An interesting theoretical question that seems more amenable to analytic study 
is the scaling behavior in the $s \downarrow 1$ limit.
That is, in the spirit of ``universality" in statistical physics, 
one can speculate that there exists an exponent $\alpha$ such that 
\[  \Psi(s) \asymp (s-1)^{-\alpha} \mbox{ as } s \downarrow 1 \]
where the value of $\alpha$ does not depend on any detailed assumptions in the model 
(worst-case or average-case; the Steiner network case or the graph network case)
but instead depends only on the the fact the we are studying the length-stretch trade-off in two-dimensional space.  
Our results  imply crude bounds on $\alpha$: 
an upper bound of $\frac{3}{4}$ for $\Psave$ (Corollary \ref{C34}) 
and  $\frac{5}{4}$ for $\Psworst$ (\ref{54bound}), 
and a lower bound of $\frac{3}{8}$ for $\Psave$ and hence for $\Psworst$ also (Proposition \ref{P38}).  
But considering these $s \to 0$ limits is considering increasingly dense networks (``covering the countryside in tarmac", we say in talks), 
  which is hard to motivate.

\subsection{What is already known?} 
\label{sec:already}
The literature on geometric spanner networks focusses on worst-case bounds on stretch (and many other statistics of networks) produced by algorithmic procedures from arbitrary configurations, works in the ``graph network" setting rather than our Steiner network setting, and emphasizes hop length more than route length.
So while the techniques of that field are clearly relevant, it is not so easy to directly apply their results to the study of $\Psave$ and $\Psworst$.

One relevant result in 
\cite{MR2289615} is  Theorem 15.2.16, which says that for small $s>0$ one can construct 
$(1+s)$-spanners such that (amongst other properties) 
the network length is bounded by 
$O(s^{-4} \times \mbox{(length of MST)})$, 
where MST denotes the minimum spanning tree on the given configuration.
It is well-known and elementary 
(see e.g. \cite{MR1422018} section 2.2)
that in the worst case the length of MST is $O(n)$, so the theorem mentioned above implies
\begin{equation}
\Psworst(s) = O((s-1)^{-4}) \mbox{ as } s \downarrow 1
\label{C1}
\end{equation}
and in particular that 
$\Psworst(s) < \infty$ for all $s > 1$.

A remarkable result 
\cite{MR1134449}
is that the Delaunay triangulation is always a $t$-spanner for 
$t = \frac{2\pi}{3 \cos \pi/6} \approx 2.42$.
Because the length of a Delaunay triangulation is {\em not} $O(n)$ in the worst case, 
this does not directly help us bound $\Psworst$.
But in the random model, a classical result (\cite{MR0279853} page 113) shows the limit normalized length of the 
 Delaunay triangulation equals
$\frac{32}{3\pi} = 3.40...$.
So we get a numerical bound 
\begin{equation}
\Psave(2.42) \leq 3.40.  
\label{C2}
\end{equation}

More generally, there are known bounds \cite{bose2013} on stretch for the 
well-studied one-parameter {\em $\theta$-graph} family of graph networks.
As above, their lengths are  not $O(n)$ in the worst case, so
this does not directly help us bound $\Psworst$.  
But again we can calculate mean lengths in the random model, and so deduce 
explicit upper bounds on $\Psave(s_m)$ for a certain sequence 
$s_m \downarrow 1$ (section \ref{sec-5}). 
However for small $s$ we get better bounds, in our Steiner network setting, from the construction in section \ref{sec-99}.

 \subsection{Other statistics for route-length efficiency}
In defining a statistic to summarize the effectiveness of a network in providing short routes,  
one may be more interested in the typical value of
\[ 
R(v,w) = 
\frac{r(v,w)}
{d(v,w)}
 - 1  \] 
than in the maximum value used in the definition of {\em stretch}. 
One might first consider the summary statistic 
$\ave_{v,w} R(v,w)$, which somewhat counter-intuitively can easily be made 
very small for large $n$ \cite{MR2411811,abs-1305-4170}.
It is argued in \cite{MR2791668} that the most appropriate summary statistic $R$ is  
defined as follows. 
For each distance $d$, set 
$\rho(d) = $  average of $R(v,w)$ over city-pairs at distance approximately $d$; 
then let $R$ be the maximum of $\rho(d)$ as $d$ varies.  
The trade-off between $R$ and normalized length, in the average-case setting, is discussed in \cite{MR2791668}, and the motivation for
the present paper was to make a connection with the topic of spanner networks.

\section{Existence of the limit functions $\Psi$}
\label{sec-2}
In this section we use a subadditivity argument in the spirit of
\cite{MR1422018,MR1632875}
to prove 
existence of the limit functions $\Psi$.
Note that in the most familiar kind of spatial subadditivity argument 
a big square is divided into small subsquares, and optimal solutions
on subsquares are used to construct some near-optimal solution on the big square.
We argue in the opposite direction: use an optimal solution on the big square to construct 
near-optimal solutions on subsquares. 
This leads to the ``superadditive" inequalities (\ref{ank}, \ref{bnk}).

Fix $1 < s < \infty$ and let $a_n$ be
the worst-case value 
(over configurations $\bz = (z_1,\ldots,z_n)$ of cities 
in the square of area $n$) 
of the length of the shortest network on $\bz$ with stretch $\leq s$.
We shall prove existence of the limit
\begin{equation}
\Psworst(s) := \lim_n n^{-1}a_n  \leq \infty . 
\label{Psworst-def}
\end{equation}

We will first argue
\begin{equation}
\sfrac{a_n}{n} \leq \sfrac{a_{nk^2}}{nk^2} + \sfrac{4}{\sqrt{n}}, \quad 
n \geq 1, \ k \geq 2 .
\label{ank}
\end{equation}
Fix $n$ and $k$.
Let $\bz$ be a configuration in the area-$n$ square attaining $a_n$.
Take $k^2$ copies of this configuration, and translate each to construct a configuration $\bz^*$ of $nk^2$ 
cities in the square of area $nk^2$.  
By definition of $a_{nk^2}$ there is a network on $\bz^*$ with stretch $ \leq s$ and with length $\leq a_{nk^2}$.
Add to this network the four boundary edges of each of the $k^2$ subsquares
(so we get two copies of each edge interior to the big square).
We now have a network $\NN^*$ whose length $c_{n,k}$ satisfies 
$c_{n,k} \leq a_{nk^2} + 4 n^{1/2} k^2$. 
Consider the restriction of this network to one of the subsquares.
The length of the restricted network may depend on the subsquare, but there must be at least one subsquare 
$Q$ 
for which this length of the restricted network $\NN_Q$ is at most the average $c_{n,k}/k^2$.
Routes in the network $\NN^*$ between cities of 
$Q$ might go outside $Q$, but replacing these external segments by the boundary edges of $Q$ can only shorten the route length,
so $\NN_Q$ has stretch at most $s$. 
But $\NN_Q$ defines (by translation) a network on the original configuration $\bz$, and so 
$a_n \leq c_{n,k}/k^2$, which gives (\ref{ank}).

Deducing existence of a limit from (\ref{ank}) is one of many variants of routine ``subadditivity" arguments,
as follows.
First note that $(a_n)$ is increasing; indeed 
\[
a_{n+1} \geq a_n \sqrt{\frac{n+1}{n}} 
\]
by adding an arbitrary city to the configuration attaining $a_n$ and rescaling.
(Note this is one of many minor ways in which Steiner networks are technically more tractable than graph networks).
Next define 
\[ \gamma = \liminf_n a_n/n  \leq \infty \]
and use monotonicity to show
\[ \gamma = \liminf_k \frac{a_{nk^2}}{nk^2}, \quad \mbox{ for each fixed } n . 
\]
Then (\ref{ank}) shows 
\[ a_n/n \leq \gamma + 4n^{-1/2} \]
and so $\limsup_n a_n/n \leq \gamma$, meaning that indeed 
$\lim_n a_n/n = \gamma$.

This argument shows $\Psworst(s) \leq \infty$ exists.  
As mentioned at (\ref{C1}), existing results then imply 
$\Psworst(s) < \infty$ for all $s>1$; this alternatively could be derived from 
the more elementary constructions in our section \ref{sec-4}.

For the random model we use the same construction with a Poissonized number of random points.
Fix $s$ again.
Let $b_n$ be the expectation of the length of 
the shortest network with stretch $\leq s $
over $n$ uniform random cities in the {\em unit} square.
So $n^{1/2} b_n$ is the corresponding expectation in the area-$n$ square.
We shall prove existence of the limit
\begin{equation}
\Psave(s) := \lim_n 
\frac{n^{1/2}b_n}{n} .
\label{Psave-def}
\end{equation}
Write $N(t)$ for a random variable with Poisson($t$) distribution and write 
\[ \beta_t = \Ex b_{N(t)} . \]
Take a Poisson point process (rate $1$ per unit area) of cities on the whole plane.
Now $t^{1/2} \beta_t$ is the expectatation of the length of
the shortest network with stretch $\leq s $
on the Poisson cities in an area-$t$ square.
Consider partitioning a square of area $tk^2$ into $k^2$ subsquares of area $t$.
Repeating the argument for (\ref{ank}),
now using a random subsquare, 
gives an inequality analogous to (\ref{ank}):
\begin{equation}
\frac{\beta_t}{t^{1/2}} \leq \frac{\beta_{tk^2}}{t^{1/2}k} + \frac{4}{t^{1/2}}, \quad 
0<t<\infty, \ k \geq 2 . \label{bnk}
\end{equation}
Using the fact that $\beta_t$ is increasing in $t$, we can repeat the 
``subadditivity" argument to show existence of the limit 
\begin{equation}
\lim_t t^{-1/2} \beta_t = \gamma^* \leq \infty . 
\end{equation}
The ``average case better than worst case" inequality 
$n^{1/2}b_n \leq a_n$ 
and the fact $\Psworst(s) < \infty$ 
now easily imply $\gamma^* < \infty$.

To finish we need a routine ``dePoissonization" argument to show that 
(\ref{bnk}) and monotonicity of $b_n$ imply 
\begin{equation}
\lim_n n^{-1/2} b_n = \gamma^* .
\label{Bnk}
\end{equation}
First fix $\eps > 0$ and consider $
t_n/n \to 1+ \eps$.
Then
\[ \beta_{t_n} \geq b_n \Pr(N(t_n) \geq n) = b_n (1 - o(1)) \]
so 
\[ \limsup_n n^{-1/2} b_n \leq \limsup_n n^{-1/2} \beta_{t_n}
= (1+\eps)^{1/2} \gamma^* \]
and the upper bound for (\ref{Bnk}) follows.
Next, the following property of the Poisson distribution
\[ \max_{i \geq n} \frac{\Pr(N((1-\eps)n)=i)} {\Pr(N(n)=i)} \to 0 \]
implies
\[ \Ex b_{N((1-\eps)n)} \ind (N((1-\eps)n) \geq n) 
= o(\Ex b_{N(n)})  = o(n^{1/2}) \]
and so
\[ \beta_{(1-\eps)n} \leq b_n + o(n^{1/2}) \]
implying
\[ \liminf_n n^{-1/2}b_n \geq \liminf_n n^{-1/2} \beta_{(1-\eps)n} 
= (1-\eps)^{1/2} \gamma^* \]
and the lower bound for (\ref{Bnk}) follows.

\subsection{Poisson process interpretation of $\Psave$}
\label{sec-Pois-interpret} 
The argument above interprets $\Psave(s)$ as an $n \to \infty$ limit of the random $n$-city model.  
By standard weak convergence arguments which we will not give here 
(see e.g. \cite{MR2235175} section 3.5 for more details in a somewhat similar setting) we can give an ``exact" interpretation 
of $\Psave(s)$ in terms of a Poisson (rate $1$) point process of cities on the infinite two-dimensional plane.  
Consider a network $\NN_\infty$ on such cities whose distribution $\mu$ is translation invariant and ergodic.  
Associated with $\mu$ are  two numbers: 
the stretch, say $S(\mu)$, and the normalized length (mean length-per-unit area), say $L(\mu)$, which is well-defined by 
translation invariance (of course these numbers might be $+ \infty$).  
Then
\begin{equation}
\Psave(s) = \inf \{L(\mu); \ \mu \mbox{ is translation invariant, } S(\mu) \leq s \} .
\label{Poisson-interpret}
\end{equation}

\section{Short networks and the Steiner constants} 
\label{sec-3}
Write $\bz_n = (z_1,\ldots,z_n)$ 
for a configuration of city positions in the square of area $n$.
Write $\ST(\bz_n)$ for the Steiner tree 
(i.e. minimum length connected network) on $\bz_n$, 
and for any network $\NN$ write 
$\len(\NN)$ for its total length.
By an easy ``superadditive" argument similar to that in section \ref{sec-2},
there exists a limit constant for worst-case normalized Steiner tree length:
\begin{equation}
\cworst:= 
   \lim_n \sup_{\bz_n} n^{-1} \len(\ST(\bz_n))  .
\end{equation}
It is known \cite{MR627537} that $\cworst \leq 0.995$ and that
(by considering the hexagonal lattice)
$\cworst \geq (3/4)^{1/4} = 0.9306$.  
 Clearly we must have 
$\Psworst(s) \geq \cworst$ for all $s$, and this inequality must persist in the limit 
(which exists by monotonicity):
$\lim_{s \to \infty} \Psworst(s) \geq \cworst$.  

Turning to the
average-case setting, 
it follows from the general theory of subadditive Euclidean functionals 
\cite{MR1422018,MR1632875} that there exists a limit constant $\cave$ such that
 \begin{equation}
 n^{-1} \len(\ST(Z_1,Z_2,\ldots,Z_n)) \to \cave \mbox{ in } L^1  
 \end{equation} 
where the $(Z_i)$ are independent uniform random in the area-$n$ square.
 As above, we clearly have 
 $\lim_{s \to \infty} \Psave(s) \geq \cave$.  
 It is natural to guess (but not obvious) that these limit inequalities are really equalities.
 This guess is correct, as an immediate corollary of the following estimate for arbitrary city configurations. 
 \begin{Proposition}
 \label{P0}
 There exists a  function $\delta(s) \leq \infty $  with $\lim_{s \to \infty} \delta(s) = 0$,
 and a function $K(s) < \infty$,  
 such that  for all $1<s<\infty$, all $n \geq K(s)$ and all 
 city configurations $\bz_n$ in the area-$n$ square,
 there exists a network $\NN$ connecting cities $\bz_n$ such that
\[ \stretch(\NN) \leq s; \quad
 n^{-1} (\len (\NN) - \len( \ST(\bz_n))) \leq \delta(s) .\]
 \end{Proposition}
 \begin{Corollary}
 \[ \mbox{$\lim_{s \to \infty} \Psworst(s) = \cworst$ and  $\lim_{s \to \infty} \Psave(s) = \cave$.} \]
 \end{Corollary} 
 The idea of the proof is to partition the area-$n$ square into rectangles containing at most $K$ cities, and then use a crude construction (Lemma \ref{L2}) of networks on $K$ cities.  
We will set up some notation, state the lemma, give the reduction of the Proposition to the lemma, and then prove the lemma.

Fix $K \geq 0$. 
Let $A$ be a rectangle; write $\partial A$ for its boundary, so that $\len (\partial A)$ is its boundary length.  
Let $y_1, \ldots, y_K$ be an arbitrary configuration of $K$ cities in $A$.  
Consider a network $\NN = \NN(A)$ in $A$ which includes the boundary $\partial A$ and links the cities to the boundary.
For such a network define
\begin{equation}
 \stretch^*(\NN) = \max_{y \neq y^\prime} \sfrac{\mbox{route-length from $y$ to $y^\prime$}}{d(y,y^\prime)}  \label{R-bar} 
\end{equation}
where $y$ and $y^\prime$ run over the cities and over points of $\partial A$.
\begin{Lemma}
\label{L2}
Let $\hat{\bt}$ be the Steiner tree on the cities $y_1, \ldots, y_K$ in a rectangle $A$ and (possibly) other cities outside $A$.
Let $\bt$ be the intersection of $\hat{\bt}$ with $A$.
There there exists a network $\NN$ in $A$ containing $\bt$ 
and  $\partial A$ and linking the cities to $\partial A$,
such that 
\[ \stretch^*(\NN) \leq \rho(K); \quad 
\len(\NN) - \len(\bt)  \leq 2 \len(\partial A) \] 
where $\rho(K) < \infty$ depends only on $K \ge 0$.
\end{Lemma}

\paragraph{Proof of Proposition \ref{P0}.} 
Fix $K$ and $n>K$.
We use a simple decomposition, the {\em multidimensional search tree} or {\em $k-d$ tree} \cite{MR805539}.  
Split the square $[0,n^{1/2}]^2$ into two rectangles using a vertical line through the city with median $x$-coordinate (if $n$ is odd) or a vertical line separating the two median $x$-coordinate cities (if $n$ is even).  
In either case, each rectangle has at most $n/2$ cities in its interior.  Separately for each rectangle, split it into two rectangles using horizontal lines through the median $y$-coordinate(s).  Now (end of stage $1$) 
we have $4$ rectangles, each with at most $n/4$ cities in its interior.   Continue recursively for $L$ stages, where $L$ is the smallest integer such that $n 4^{-L} \leq K$, to get a partition into $4^L$ rectangles, each with at most $K$ cities in its interior.  
Write $A$ for a generic rectangle in this partition.

Given a configuration $\bz_n$ in $[0,n^{1/2}]^2$, apply Lemma \ref{L2} 
(where $\hat{\bt}$ is the Steiner tree on $\bz_n$) 
 to each $A$ and the cities inside $A$ to obtain a network $\NN(A)$ satisfying
 \begin{equation}
  \stretch^*(\NN(A)) \leq \rho(K); \quad 
\len(\NN(A)) - \len(\ST(\bz_n) \cap A)  \leq 2 \len(\partial A) .
\label{NAK}
\end{equation}
 Then consider the network $\NN$ on the cities $\bz_n$  obtained as the union of networks $\NN(A)$.
Note that  the bound $\rho(K)$ on $\stretch^*(\NN(A))$ does not depend on $A$. 
For any pair of cities $z_i, z_j$, we can define a route in $\NN$ between them by 
considering the points $v_1,v_2, v_3,\ldots$ at which a straight line between them intersects boundaries of successive rectangles $A_1, A_2, A_3,\ldots$, and within each such rectangle $A$ use the shortest route in $\NN(A)$ between these boundary 
points (or the cities $z_i, z_j$ themselves, at the ends).
It follows that  $\stretch(\NN) \leq \max_A \stretch^*(\NN(A)) \leq \rho(K)$. 
Note that the intermediate rectangles may contain no cities of $\bz_n$, explaining why we must allow 
$K = 0$ in Lemma \ref{L2}.

As a preliminary to bounding $\len(\NN)$, we need to consider the total length of lines used in the original decomposition.
Include a stage $0$ in which the edges of the external boundary $\partial_0$ of 
$[0,n^{1/2}]^2$ are added.  At stage $1$, the length of lines added equals $2n^{1/2}$, and inductively at stage $j$  the length of lines added equals $2^jn^{1/2}$. 
Because each segment of these added lines (except the external boundary) is part of the boundary of exactly $2$ of the final rectangles $A$, 
\[ \sum_A \len(\partial A) = \len(\partial_0) + 2\sum_{j=1}^L 2^j n^{1/2}
= 2 n^{1/2}(2 + \sum_{j=1}^L 2^j )
= 2^{L+2} n^{1/2}. \] 
By definition of $L$ we have $n 4^{-(L-1)} > K$, giving $2^L \leq 2 n^{1/2}K^{-1/2}$, and so 
\begin{equation}
  \sum_A \len(\partial A) \leq 8 n K^{-1/2} . 
 \label{ndA}
 \end{equation}
 So
 \begin{eqnarray*}
 \len (\NN) &\leq& \sum_A \len(\NN(A)) \\
 &\leq& \sum_A \left( \len (\ST(\bz_n) \cap A) + 2 \len(\partial A) \right) \mbox{ by (\ref{NAK})}\\
 &=& \len (\ST(\bz_n)) + 2   \sum_A \len(\partial A) .
 \end{eqnarray*}
 Combining with (\ref{ndA}), 
\[ 
 n^{-1} (\len (\NN) - \len( \ST(\bz_n))) \leq 16 K^{-1/2} . \] 
We may assume $\rho(K) \uparrow \infty$ as $K \uparrow \infty$, and now
Proposition \ref{P0}  holds
for $K(s):= \max\{ K: \rho(K) \leq s\} $ 
and $\delta(s):= 16 K^{-1/2}(s)$.

\paragraph{Proof of Lemma \ref{L2}.}  
We may suppose $A$ is an $a_1 \times a_2$ rectangle, where 
$a_1 \leq a_2$.
The network $\NN$ will consist of \\
(i) $\bt$ (the intersection of $\hat{\bt}$ with $A$) \\
(ii) the boundary $\partial A$ of $A$ \\
(iii) extra edges, of total length at most $\len (\partial A)$.

\noindent
Set $m = \lfloor a_2/a_1 \rfloor$ 
and partition $A$ into $m+1$ similar $a_1 \times \frac{a_2}{m+1}$ rectangles by using 
$m$ equally spaced roads of length $a_1$. 
So the total length of these added roads is 
$ma_1 \leq a_2 \leq \frac{1}{2} \len (\partial A)$. 
So the network $\NN_0$ consisting of $\bt$ and $\partial A$ and these extra roads 
has 
$\len(\NN_0) - \len(\bt) \leq \frac{3}{2} \len(\partial A)$.
It is easy to check that this network $\NN_0$ (without using the edges of $\bt$) satisfies 
\[ \max_{y \neq y^\prime \in \partial A} \sfrac{\mbox{route-length from $y$ to $y^\prime$}}{d(y,y^\prime)} \leq 2 . \]
In particular, the $K = 0$ case of Lemma \ref{L2} holds with $\rho(0) = 2$.

Now consider the case $K \geq 1$.
To cover the possibility that $\hat{\bt}$ and hence $\bt$ is entirely in the interior of one of the subrectangles of $A$, 
add to $\NN_0$ a road to the boundary from the city closest to the boundary.
This road has length at most $\sfrac{1}{2}a_1 \leq \sfrac{1}{8} \len (\partial A)$,
and the resulting network $\NN_1$ has length 
$\len(\NN_1) - \len(\bt) \leq \frac{13}{8} \len(\partial A)$.

Now set 
\[ \eta:= 
\frac{\sfrac{3}{8} \len (\partial A)}
{K + {K \choose 2}} 
. \]
     Let $\NN$ be the network $\NN_1$ augmented as follows:
     for each
     city within distance $\eta$
      from the boundary,
     add a road from the city to the closest boundary point;
     for each pair of cities within distance $\eta$
      of each other,
     add a road directly linking them.
From the definition of $\eta$, the extra length added in this stage is at most 
$\sfrac{3}{8} \len (\partial A)$, 
and so $\NN$ satisfies the length requirement 
\[ \len(\NN) - \len(\bt)  \leq 2 \len(\partial A) \] 
in Lemma \ref{L2}.

It remain to bound $\stretch^*(\NN)$.  We quote a simple bound on Steiner tree length
(given for squares in \cite{MR2411811} Lemma 10; the extension to rectangles is straightforward).
\begin{Lemma}
Under the assumptions of Lemma \ref{L2},
\[ \len(\bt) \leq C_1(K) \len (\partial A) \]
where $C_1(K)$ depends only on $K$.
\end{Lemma}
So 
$\len(\NN) \leq (2 + C_1(K)) \len (\partial A)$ and then 
\begin{equation}
\frac{\len(\NN)}{\eta} \leq \sfrac{8}{3} (2 + C_1(K)) (K + {K \choose 2}) .
\label{222}
\end{equation}

To bound $\stretch^*(\NN)$ we need to treat several cases 
for the pairs $(y,y^\prime)$ in (\ref{R-bar}).
We have already obtained an upper bound of $2$ for the case where both points 
are on the boundary.
If the two points are at distance $\geq \eta$ apart then,
because route length $r(y,y^\prime)$ is at most network length,
$\sfrac{r(y,y^\prime)}{d(y,y^\prime)}
\leq \len(\NN)/\eta $.
If the two points are cities within distance $\eta$ then 
$\sfrac{r(y,y^\prime)}{d(y,y^\prime)} = 1$.
The only remaining case is a city $y$ within distance $\eta$ from the boundary, and a boundary point
$y^\prime$ within distance $\eta$ from $y$.  In this case, by using the edge from $y$ to the closest 
boundary point and then following the boundary we find (the worst case is near a corner)
$\sfrac{r(y,y^\prime)}{d(y,y^\prime)} \leq 3$.  
So
\[ \stretch^*(\NN) \leq \max(3, \sfrac{\len(\NN)}{\eta} ) \] 
and by (\ref{222}) we have proved Lemma \ref{L2}.

\section{Upper bounds on $\Psave$ from worst-case stretch for $theta$-graphs}
\label{sec-5}
In this section we show how to derive upper bounds on $\Psave$  from known bounds on worst-case stretch for graph networks.

The $\theta_m$-graph on a configuration is defined as follows.
At each point $z_i$, consider the natural partition of the plane into
 $m$ equal-angle cones of base angle $\theta_m = 2\pi/m$ based at $z_i$;
 the boundary lines make angles $(2 \pi i/m, 0 \le i < m)$ with the $x$-axis.
Given $z_i$ and such a cone, each point $z_j$ in the cone has an orthogonal 
projection onto the bisector line of the cone, at position $z^\prime_j$ say; create an edge $(z_i,z_j)$ for the point $z_j$ in the cone such that 
$z^\prime_j$ is closest to $z_i$.
This is a now well-known construction of graphs with low stretch, and the best known explicit bounds on stretch are given
as follows in \cite{bose2013}.  These hold for $m \ge 6$.

\begin{eqnarray*}
1 + \frac{2 \sin(\theta_m/2)}{\cos(\theta_m/2) - \sin(\theta_m/2)}&\quad & m = 0 \bmod 4 \\
\frac{\cos(\theta_m/4)}{\cos(\theta_m/2) - \sin(3\theta_m/4)} & \quad & m = 1 \mbox{ or } 3 \bmod 4 \\
1 + 2 \sin(\theta_m/2) & \quad & m = 2 \bmod 4.
\end{eqnarray*}
Writing $s_m$ for these bounds, we immediately have from (\ref{Poisson-interpret})
\begin{equation}
\Psave(s_m ) \le L_m 
\label{smLm}
\end{equation}
where $L_m$ is the mean length-per-unit-area of the $\theta_m$ graph over the rate-$1$ Poisson point process on $\Reals^2$. 
Calculating $L_m$ is in principle straightforward; indeed more detailed
calculations of various statistics in the finite-$n$ random model can be found in
\cite{morin_v}, though they do not explicitly consider the statistic $L_m$.  
We use some of the notation from \cite{morin_v}.
As observed there, the calculation is in practice easier in the case of even $m$, so we treat that case.

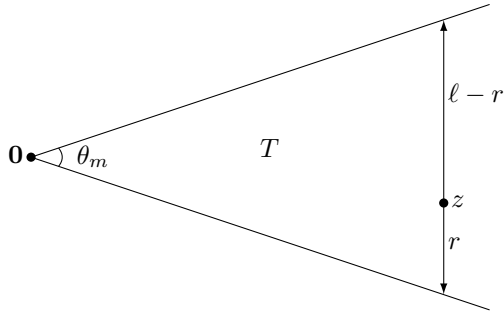
\begin{figure}
\setlength{\unitlength}{0.24in}
\begin{picture}(10,7)(0,-3)
\put(0,0){\line(3,1){10}}
\put(0,0){\line(3,-1){10}}
\put(0,0){\circle*{0.2}}
\put(9,-1){\circle*{0.2}}
\put(9,-1){\vector(0,1){4}}
\put(9,-1){\vector(0,-1){2}}
\put(-0.5,-0.12){$\origin$}
\put(9.16,-1.1){$z$}
\put(9.1,1.2){$\ell - r$}
\put(9.1,-2){$r$}
\put(5,0){$T$}
\qbezier(0.6,0.2)(0.76,0)(0.6,-0.2)
\put(1,-0.14){$\theta_m$}
\end{picture}
\caption{Parametrization of $z$ as $(r,\ell)$.}
\end{figure}

Take the origin $\origin$ as a typical point of the Poisson configuration.
An edge $(\origin,z)$ created by the defining rule applied at $\origin$ may or may not be {\em mutual}, meaning it is also  created by the rule applied 
at $z$.
We readily see the formula
\begin{equation}
L_m = m \int_C \ ||z|| 
\left(\sfrac{1}{2} \pmut(z) + \pnot(z) \right) \ dz .
\label{Lm}
\end{equation}
Here 

$C$ is a cone of base angle $\theta_m = 2\pi/m$;

$||z||$ is Euclidean distance from $z$ to the origin;

$\pmut(z)$ is the probability that (if there is a Poisson point at $z$)
there is a mutual edge  $(\origin,z)$; 
$\pnot(z)$ is the corresponding probability of a non-mutual edge created by the rule at $\origin$.  
The $\frac{1}{2}$ term avoids double-counting mutual edges.

\smallskip
\noindent
Following \cite{morin_v} we first parameterize a point $z \in C$ by a pair 
$(r,\ell) = (r(z),\ell(z))$ where, drawing the bisector horizontally in 
Fig. 2,
$r$ and $\ell - r$ are the vertical  distances from $z$ to the cone boundaries.  
If there is a Poisson point at $z$ then the rule at $\origin$ creates an edge
iff triangle $T = T(z)$ is empty of other Poisson points. 
Moreover this edge will be mutual iff a certain other triangle $T^\prime = T^\prime(z)$ is also empty.
Now \cite{morin_v}
\[ \area(T) = \alpha \ell^2; \quad \area(T^\prime \setminus T) = 
\alpha(r^2 + (\ell - r)^2) \]
where 
\[ \alpha := \sfrac{\cos (\theta_m/2)}{4 \sin(\theta_m/2)} .\] 
So
\[
\pmut(z) = \exp(- \area(T \cup T^\prime)) 
= \exp(- \alpha(\ell^2 + r^2 + (\ell - r)^2)) \]
\[ \pnot(z) = \exp( - \area(T)) - \exp(- \area(T \cup T^\prime)) 
= \exp(-\alpha \ell^2) - \exp(- \alpha(\ell^2 + r^2 + (\ell - r)^2)) 
\]
and so 
\[
\sfrac{1}{2} \pmut(z) + \pnot(z) = 
 \exp(-\alpha \ell^2) - 
\sfrac{1}{2}  \exp(- \alpha(\ell^2 + r^2 + (\ell - r)^2)) 
 . \]
Substituting into (\ref{Lm}) we now have an expression for $L_m$, for even 
$m \ge 6$.
It is straightforward to show the asymptotics 
\[ L_m = \Theta( m^{3/2}) \quad s_m - 1  = \Theta(m^{-1}) \]
and then from (\ref{smLm}) we have 
\begin{equation}
\Psave(s) = O((s-1)^{-3/2} ) \mbox{ as } s \downarrow 1
 . \label{s32}
\end{equation}
We will see in section \ref{sec-99} that we can improve this bound when using 
Steiner networks instead of graph networks.

\section{Upper bounds via a ``freeways and access roads" construction}
\label{sec-4}  
In our ``Steiner network" setting we can get explicit bounds on 
$\Psworst$ via
 elementary constructions using parallel ``freeways" in different directions, 
with ``access roads" linking cities to nearby freeways.  
We give details in the simplest setting in section \ref{sec:sq}, and state 
a more general result in section \ref{sec-Tamar}.

\subsection{Constructions based on a square grid of roads}
\label{sec:sq}
\begin{Proposition} 
\label{P123}
\begin{eqnarray}
\Psworst(2) & \leq & 4 \label{sq-1} \\
\Psworst(\sfrac{3}{2} ) & \leq & 4 \sqrt{2} \label{sq-2} \\
\Psworst(\sqrt{2} ) & \leq & 4 \sqrt{3}  .\label{sq-3} 
\end{eqnarray}
\end{Proposition}
{\bf Proof.} 
Fix $0<t_\infty<\infty$ and choose $t = t(n) \to t_\infty$ such that  $n^{1/2}/t(n)$
is an integer $m = m(n)$.  
First construct a network of {\em grid roads} which partition the region 
$[0,n^{1/2}]^2$ into $m^2$ squares of side-length $t$.   These grid roads 
(including the boundary of $[0,n^{1/2}]^2$)
have total length
\[ n^{1/2} \times 2(m + 1) \sim 2n/t_\infty . \] 
Next, for each city construct a north-south (N-S) and an east-west (E-W) road
through the city and across the square containing the city.  These {\em access
roads} have total length 
$2tn$.

We now study the network $\NN_n^1$ thus constructed. 
We have already seen that
\begin{equation}
n^{-1} \len(\NN_n^1) \to 2(t_\infty + \sfrac{1}{t_\infty}) 
\label{LN1}
\end{equation}
so we need to bound the stretch.
Note that in a right angle triangle with side-lengths $a, b$ and $ c =
\sqrt{a^2+b^2}$ 
we have 
\[ \sfrac{a+b}{c} \leq \sqrt{2}   . \]
Thus to show that a city-pair $(i,j)$ has $\frac{r(i,j)}{d(i,j)} \leq \sqrt{2} $ it is enough to 
show that (supposing w.l.o.g. that city $j$ is to the south-east of city $i$) 
there is a route from $i$ to $j$ using only southward and eastward roads. 
But, consulting Fig. 2, this is clearly true in the three cases\\ (i) the two
cities are in the same square (as $a$ and $b$) \\
(ii) the two cities are in different rows and different columns (as $a$ and $c$).\\
(iii) the two cities are in adjacent squares (as $a$ and $d$). \\
So it remains to consider the final case \\
(iv) the two cities are in squares in the same column (say) separated by some number
$k \geq 1$ of squares.\\

\begin{figure}
 \setlength{\unitlength}{0.4in}
 \begin{picture}(7,6.6)(-2,-0.2)
\put(0,0){\line(1,0){6}}
\put(0,1){\line(1,0){6}}
\put(0,2){\line(1,0){6}}
\put(0,3){\line(1,0){6}}
\put(0,4){\line(1,0){6}}
\put(0,5){\line(1,0){6}}
\put(0,6){\line(1,0){6}}
\put(0,0){\line(0,1){6}}
\put(1,0){\line(0,1){6}}
\put(2,0){\line(0,1){6}}
\put(3,0){\line(0,1){6}}
\put(4,0){\line(0,1){6}}
\put(5,0){\line(0,1){6}}
\put(6,0){\line(0,1){6}}
\put(2.8,4){\line(0,1){1}}
\put(2,4.3){\line(1,0){1}}
\put(2.8,4.3){\circle*{0.17}}
\put(2.5,4.1){b}
\put(2.2,4){\line(0,1){1}}
\put(2,4.7){\line(1,0){1}}
\put(2.2,4.7){\circle*{0.17}}
\put(2.22,4.78){a}
\put(4.4,1){\line(0,1){1}}
\put(4,1.4){\line(1,0){1}}
\put(4.4,1.4){\circle*{0.17}}
\put(4.46,1.49){c}
\put(3.6,4){\line(0,1){1}}
\put(3,4.6){\line(1,0){1}}
\put(3.6,4.6){\circle*{0.17}}
\put(3.6,4.69){d}
 \end{picture}
 \caption{All the grid roads and some of the access roads in $\NN^1_n$.}
 \end{figure}
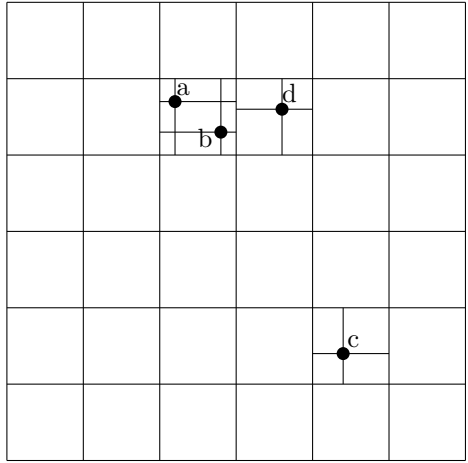
 
The remainder of the argument rests upon being able to recognize, within case (iv), which city positions $(v,w)$ maximize the ratio $r(v,w)/d(v,w)$.  
In the context of the square grid, these ``worst situations" are intuitively clear, and we will state them without proof. 
It turns out (see Fig. 3, left diagram) that the worst situation in case (iv) is where $k = 1$, this intervening square contains no cities, 
and the two cities are (arbitrarily close to) the centers of the north and the south
edges of the intervening square 
(as $e$ and $f$).  
In this situation $r(v,w)/d(v,w) = 2$, so this is an upper bound for case (iv).
 Thus the networks $\NN_n^1$ have 
$\stretch(\NN_n^1) \leq 2$. 
Consulting (\ref{LN1}), we can choose $t_\infty = 1$ so that 
 $\len(\NN_n^1) \sim 4n$,  establishing (\ref{sq-1}).

 \begin{figure}
 \setlength{\unitlength}{0.4in}
 \begin{picture}(7,7.3)(-0.4,-0.2)
\put(0,0){\line(1,0){2}}
\put(4,0){\line(1,0){2}}
\put(8,0){\line(1,0){2}}
\put(0,2){\line(1,0){2}}
\put(4,2){\line(1,0){2}}
\put(8,2){\line(1,0){2}}
\put(0,4){\line(1,0){2}}
\put(4,4){\line(1,0){2}}
\put(8,4){\line(1,0){2}}
\put(0,6){\line(1,0){2}}
\put(4,6){\line(1,0){2}}
\put(8,6){\line(1,0){2}}
\put(0,0){\line(0,1){6}}
\put(2,0){\line(0,1){6}}
\put(4,0){\line(0,1){6}}
\put(6,0){\line(0,1){6}}
\put(8,0){\line(0,1){6}}
\put(10,0){\line(0,1){6}}
\put(1.0,4.07){\circle*{0.14}}
\put(0.9,4.22){e}
\put(1.0,1.93){\circle*{0.14}}
\put(0.77,1.6){f}
\put(5.0,0){\line(0,1){6}}
\put(4,1){\line(1,0){2}}
\put(4,3){\line(1,0){2}}
\put(4,5){\line(1,0){2}}
\put(4.5,4.07){\circle*{0.14}}
\put(4.32,4.24){g}
\put(4.5,1.93){\circle*{0.14}}
\put(4.28,1.64){h}
\put(8.67,0){\line(0,1){6}}
\put(9.33,0){\line(0,1){6}}
\put(8,0.67){\line(1,0){2}}
\put(8.0,1.33){\line(1,0){2}}
\put(8,2.67){\line(1,0){2}}
\put(8.0,3.33){\line(1,0){2}}
\put(8,4.67){\line(1,0){2}}
\put(8.0,5.33){\line(1,0){2}}
\put(9.0,4.07){\circle*{0.14}}
\put(8.9,4.22){e}
\put(9.0,1.93){\circle*{0.14}}
\put(8.77,1.6){f}
\put(0.87,6.21){$\NN^1$}
\put(4.87,6.21){$\NN^2$}
\put(8.87,6.21){$\NN^3$}
 \end{picture}
 \caption{Networks with grid roads and interior roads, and the positions maximizing 
$r(v,w)/d(v,w)$.   The access roads in Fig. 2 are present but not shown; they are not helpful for these extremal positions.}
 \end{figure}
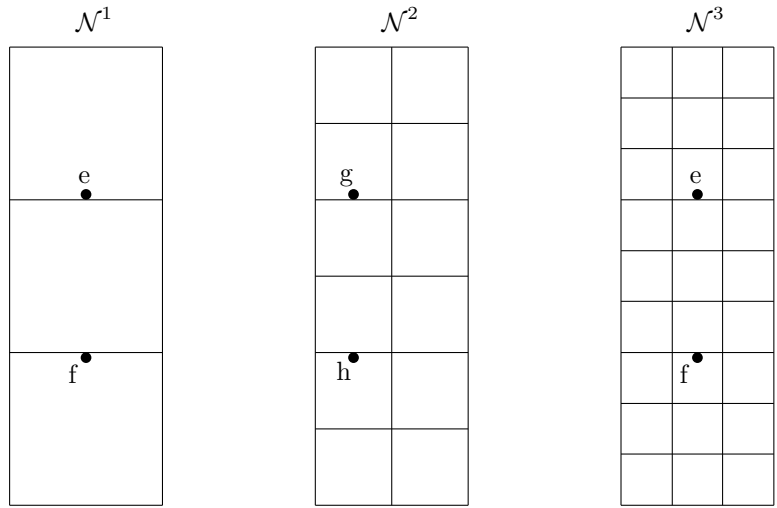

Now consider the networks $\NN_n^2$ (Fig. 3, center diagram)  obtained from $\NN_n^1$ (left diagram) by adding, for each
square, 
the N-S and the E-W {\em interior roads} across the square through the center of the square.  
Now the case (iv) worst situation is where (as $g$ and $h$ in center diagram) the two cities are
arbitrarily close to a quarter of the way along the north and the south edges of
the intervening square.  
In this situation $r(g,h)/d(g,h) = 3/2$, so this is an upper bound for case (iv).
 That is, 
$\stretch(\NN_n^2)  \leq 3/2$.  The total extra network length is $2n/t$, so  $n^{-1} \len(\NN_n^2)
\to 2(t_\infty + \sfrac{2}{t_\infty})$. 
Choosing $t_\infty = \sqrt{2}$ gives $n^{-1} \len(\NN_n^2) \to 4\sqrt{2}$ and establishes
(\ref{sq-2}).

Finally consider the networks $\NN_n^3$ obtained from $\NN_n^1$ by adding, for each
square, 
two N-S and two E-W {\em interior roads} partitioning the square into nine equal subsquares.  
Here the case (iv) worst situation is where (as $e$ and $f$ in Fig. 3, right diagram) the two cities are
arbitrarily close to half of the way along the north and the south edges of the
intervening square.  
In this situation $r(e,f)/d(e,f) = 4/3$, so this is an upper bound for case (iv).  
But here $4/3$ is less than the bound $\sqrt{2} $ from the other cases.
So  
$\stretch(\NN_n^3)   \leq \sqrt{2} $.  The total extra network length (relative to
$\NN_n^1$) is $4n/t$, so  $n^{-1} \len(\NN_n^3) \to 2(t_\infty + \sfrac{3}{t_\infty})$. 
Choosing $t_\infty = \sqrt{3}$ gives $n^{-1} \len(\NN_n^3) \to 4\sqrt{3}$ and establishes
(\ref{sq-3}).

\subsection{A generalization}
\label{sec-Tamar}
The constructions above were based on horizontal and vertical lines, distance $t$ apart. 
One can regard that as the $m = 2$ case of the line pattern
with $m$ lines through the origin at angles $\pi/m$ apart, 
each duplicated by parallel lines distance $t$ apart.  
Analogous network constructions based on this line pattern were studied in
the Master's thesis \cite{lando}, where it was shown that, for fixed $1<s<2$, the construction gives an $s$-spanner with total length bounded by 
the quantity 
$\Psi^*(s)$ below, 
which is therefore an upper bound on 
$\Psworst(s)$.
\begin{Theorem}
\label{T-Tamar}
For $1<s<2$ set 
\[ \phi_s = \sfrac{\pi}{2} - \sin^{-1}\left(\sfrac{1}{s}\right) \]
\[ \Psi^*(s) = \frac{2 \lceil \sfrac{\pi}{\phi_s} \rceil 
\sqrt{(1 + \lceil \sfrac{1}{s-1} \rceil) \tan \phi_s}
}{\sin \phi_s} . \]
Then 
$\Psworst(s) \leq \Psi^*(s)$.
\end{Theorem}
In particular, as $s \downarrow 1$ we have 
$\phi_s \sim \sqrt{2(s-1)}$ and then 
$\Psi^*(s) \sim 2^{1/4} \pi (s-1)^{-5/4}$, 
so
\begin{equation}
\Psworst(s) = O((s-1)^{-5/4}) \mbox{ as } s \downarrow 1. 
\label{54bound}
\end{equation}
We will not repeat the proof of Theorem \ref{T-Tamar} here.

\section{Upper bounds by putting a road in every cone}
\label{sec-99}
We first show (Proposition \ref{P-theta-dense}) that one can achieve a given stretch $s>1$ by insisting that the network has the property of containing
roads from each city within each cone of appropriate base angle $\theta_s$. 
In section \ref{sec-ub2} we show how, in the random model, it is easy to construct networks with the desired property
whose expected length can be calculated; this leads to bounds on 
$\Psave(s)$, stated in Proposition \ref{P7}.
This idea is quite similar to the notions of $\theta$-graph 
from section \ref{sec-5} and of {\em Yao graph} \cite{yao1982}.
But by using Steiner networks instead of graph networks we obtain in 
 Proposition \ref{P7} a bound which (for small $s-1$, at any rate) improves the bound (\ref{s32}) derived from $\theta$-graphs.

To spotlight the essential difference between graph networks and Steiner networks here, 
Fig. 5 (copied from  Fig. 5 of \cite{bose_optimal}) illustrates a worst-case configuration for stretch of $\theta$-graphs: a route from $w$ to $u$ must go via $v$ or $v^\prime$.   
In our construction, there would be a line from $w$ which meets the line 
$(v,u)$ somewhere near $u$.
 It seems plausible that one can get bounds on $\Psworst(s)$ in a similar way, 
adapting other methods from \cite{MR2289615}, and perhaps improve on
Theorem \ref{T-Tamar},
but we have not investigated this question carefully.

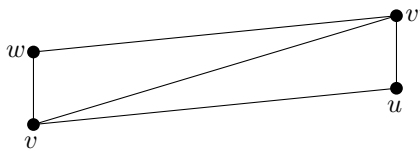
\begin{figure}
 \setlength{\unitlength}{0.095in}
 \begin{picture}(52,15)(-16,-3)
\put(0,0){\circle*{0.7}}
\put(20,2){\circle*{0.7}}
\put(0,4){\circle*{0.7}}
\put(20,6){\circle*{0.7}}
\put(0,0){\line(0,1){4}}
\put(20,2){\line(0,1){4}}
\put(0,0){\line(20,2){20}}
\put(0,4){\line(20,2){20}}
\put(0,0){\line(20,6){20}}
\put(-0.5,-1.3){$v$}
\put(-1.5,3.7){$w$}
\put(19.5,0.7){$u$}
\put(20.5,5.7){$v^\prime$}
\end{picture}
\caption{A bad configuration for a $\theta$-graph.}
\end{figure}

\subsection{The construction}
Given a point $z$ in the plane and angles 
(relative to $x$-axis, as usual)
$\phi$ and $\theta$, 
write 
$\cone(z,\phi,\phi + \theta)$ 
for the cone bounded by the two rays from $z$ at angles $\phi$ and $\phi + \theta \bmod 2 \pi$.
Fix $0 < \theta < \pi/2$.
Consider a {\em graph} network on a given configuration of cities.
Call such a network {\em $\theta$-dense} if for each city $z$ 
and each $\phi$, if there exists another city in 
$\cone(z,\phi,\phi + \theta)$,
then there exists a road from $z$ to some city in that cone.  
One can find analogs of Proposition \ref{P-theta-dense} below for finite configurations, 
but it is simpler 
(and sufficient for our purposes) 
to work under the assumption 
\begin{equation}
\mbox{for each city $z$ and each $\phi$, the $\cone(z,\phi,\phi + \theta)$ contains another city}
\label{local-finite}
\end{equation}
which of course cannot hold for any finite configuration but does hold for the Poisson process on the infinite plane.
\begin{Proposition}
\label{P-theta-dense}
Consider a locally finite configuration on the plane satisfying (\ref{local-finite}), and consider a
$\theta$-dense graph network on the configuration.  Then its stretch 
(considered as a Steiner network)
is at most $\frac{1}{\cos \frac{\theta}{2}}$.
\end{Proposition}
{\bf Proof.} 
Fix two cities, w.l.o.g. at 
$(0,0)$ and $(x_0,0)$, where $x_0 > 0$.
We first show that the Proposition can be reduced to the following lemma.
\begin{Lemma}
\label{L-phi}
Under the hypotheses of Proposition \ref{P-theta-dense},
there exist $- \theta < \phi < 0$ and 
 a route from $(0,0)$ to $(x_0,0)$ 
such that the angle of each segment lies in the range $[\phi, \theta + \phi]$.
\end{Lemma}
Each segment of the Lemma \ref{L-phi} route can be visualized as one edge of a triangle whose two other ``virtual edges" are at angles 
$\phi$ and $\theta+\phi$. 
So the length of the route is upper bounded by the length of the path of virtual edges 
for each such triangle 
(this is a path in the plane, not a route in the network).
The length of this path is the length of the path in the plane which goes
from $(0,0)$ to $(x_0,0)$ by using a line of angle $\phi$ followed by a line of angle $\theta + \phi$.
The length of this path is maximized (as $\phi$ varies) when $\phi = - \theta/2$ in which case the length 
equals $x_0/\cos(\theta/2)$,
establishing the Proposition \ref{P-theta-dense}  bound on stretch.

{\bf Proof of Lemma \ref{L-phi}.}
Fix $\phi \in (-\theta, 0)$.
Define a {\em lower route} from the city $(0,0)$ to some point on the line 
$\{(x_0,y): - \infty < y < \infty\}$
via the simple procedure: 
$v_0 = (0,0)$, and inductively

from $v_i$, follow the road to $v_{i+1}$, where $v_{i+1}$ is chosen so that the angle of the segment 
$(v_i,v_{i+1})$ is the lowest possible value in $[\phi,\phi+\theta]$ amongst all roads from $v_i$.

At each step there is some possible choice, by assumption (\ref{local-finite}) and the assumption of $\theta$-dense.   Stop the route where it crosses the line  $\{(x_0,y): - \infty < y < \infty\}$.

Define the analogous {\em upper route} using the maximum possible angle at each step.
It is easy to check that the upper route lies (weakly) above the lower route.
In particular, the routes are stopped at two points 
$(x_0,y^R_{\mbox{\tiny lower}})$ and 
$(x_0,y^R_{\mbox{\tiny upper}})$ 
where 
$y^R_{\mbox{\tiny lower}} \leq y^R_{\mbox{\tiny upper}}$.
These are {\em eastward} routes, but we can also define the analogous {\em westward} routes,
which start at city $(x_0,0)$ and are stopped at points 
$(0,y^L_{\mbox{\tiny lower}})$ 
and 
$(0,y^L_{\mbox{\tiny upper}})$ 
where 
$y^L_{\mbox{\tiny lower}} \leq y^L_{\mbox{\tiny upper}}$.
The roads in these routes 
are constrained to have angles in the same interval 
$[\phi,\phi+\theta]$ 
as in the eastward routes,

To establish the lemma it is enough to show
\begin{equation}
\mbox{one of the eastward routes meets one of the westward routes at some point}
\label{**}
\end{equation}
because then the route from $(0,0)$ to $(x_0,0)$ 
(switching between eastward and westward routes at the meeting point)  satisfies the conclusion of the lemma.
Clearly (\ref{**}) holds in the following cases: \\
(i) $0 \in [y^R_{\mbox{\tiny lower}}, y^R_{\mbox{\tiny upper}}]$ \\
(ii) $0 \in [y^L_{\mbox{\tiny lower}}, y^L_{\mbox{\tiny upper}}]$ \\
(iii) $y^L_{\mbox{\tiny upper}} < 0$ and $ y^R_{\mbox{\tiny upper}} < 0$ \\
(iv) $y^L_{\mbox{\tiny lower}} < 0$ and $y^R_{\mbox{\tiny lower}} < 0$.

\noindent
There remain two symmetric cases; w.l.o.g. we take the case\\
(v) $ y^R_{\mbox{\tiny upper}} < 0$ and $y^L_{\mbox{\tiny lower}} > 0$.\\
The argument so far uses a fixed $\phi$; now we exploit our freedom to choose $\phi$.
Rewrite $y^L_{\mbox{\tiny lower}}$ as $y^L(\phi)$ and 
rewrite $ y^R_{\mbox{\tiny upper}}$ as $y^R(\phi)$.
We are working in the case:
there exists $\phi_0 \in (- \theta, 0)$ such that 
$y^L(\phi_0) > 0$ and $y^R(\phi_0) < 0$.

Consider the eastward lower route for a given $\phi$. 
The route has some lowest angle, say $\hat{\phi} \geq \phi$.
As $\phi$ increases, the route does not change 
(and so $y^R(\phi)$ does not change)
until $\phi$ reaches $\hat{\phi}$,
at which stage $y^R(\cdot)$ may change but can only increase.

By considering $\phi$ arbitrarily close to $0$, either there is a road 
from $(0,0)$ to $(x_0,0)$
(in which case the result is trivial)
or else $y^R(\phi) > 0$.
It follows that there exists some 
$\phi^* \in [\phi_0,0)$ such that
$y^R(\phi^*) \leq 0$ but 
$y^R(\phi^* + \eps) \geq 0$ for all sufficiently small $\eps > 0$.
Now consider the two eastward lower routes for $\phi^*$ and for $\phi^* + \eps$. 
The westward upper route for $\phi^*$ must meet one of those eastward routes, so
the conclusion of the lemma holds for $\phi^*$.

\subsection{An upper bound on $\Psave(s)$}
\label{sec-ub2}
As in section \ref{sec-Pois-interpret}  we work with the Poisson process of cities on the infinite plane.
There are several ways one might try to use Proposition \ref{P-theta-dense}; we will just treat one of the simplest.
Fix $k \geq 2$.
For $0 \leq i \leq k-1$ define a network $\NN_i$ by:
\begin{quote}
for each city $z$, create a road as a line segment from $z$ to its closest neighbor city in 
$\cone(z,i \pi /k, (i+1)\pi/k)$, and 
a road to its closest neighbor city in 
$\cone(z,\pi + i \pi /k, \pi + (i+1)\pi/k)$.
\end{quote}
Network $\NN_i$ has a certain normalized length (mean length per unit area) $L_k$, which by rotational symmetry  of the Poisson point process does
not depend on $i$.   A calculation below will show
\begin{Lemma}
\label{LLk}
\begin{equation}
L_k =      \sqrt{2k} - \sfrac{1}{4} \pi^{1/2} \int_0^{\pi/k}  
[ \sfrac{\pi}{k} - \cos \omega \sin \omega + \sfrac{\sin^2 \omega}{\tan \pi/k}
]^{-3/2}  \ d \omega . 
\label{Lkformula}
\end{equation}
\end{Lemma}
Construct a network $\NN$ as the union of $\NN_i$ over $0 \leq i \leq k-1$.
Its normalized length equals $kL_k$.
And it is clearly $\theta$-dense for $\theta = 2 \pi/k$, so Proposition \ref{P-theta-dense} bounds its stretch by $1/\cos (\pi/k)$.  
In other words, using (\ref{Poisson-interpret})
\begin{Proposition}
\label{P7}
For each $k \geq 2$, 
$\Psave( \sfrac{1}{\cos \frac{\pi}{k}} ) \leq k L_k$.
\end{Proposition}
In particular, (\ref{Lkformula}) shows that $kL_k \leq 2^{1/2} k^{3/2}$, 
whereas $ \sfrac{1}{\cos \frac{\pi}{k}} - 1 \sim \sfrac{\pi^2}{2k^2}$ as $k \to \infty$, 
implying
\begin{Corollary}
\label{C34}
$\Psave(s) \leq (2^{-1/4} \pi^{3/2}  + o(1)) (s-1)^{-3/4}$ as $s \downarrow 1$.
\end{Corollary}
{\bf Proof of Lemma \ref{LLk}.} 
Write $\mathbf{0}$ for the origin.
Consider a position measured in polar coordinates as $(r,\omega)$ with $0 < \omega < \pi/k$. 
So
$(r,\omega) \in \cone(\mathbf{0},0,\pi/k)$
and
$\mathbf{0} \in \cone((r,\omega), \pi, \pi + \pi/k)$. 
Suppose there are cities at $\mathbf{0}$ and at $(r,\omega)$, with other cities at the points of a Poisson process.  Define
\begin{eqnarray*}
p(r,\omega) &=& \Pr((r,\omega) \mbox{ is nearest city to $\mathbf{0}$ in $\cone(\mathbf{0},0,\pi/k)$})\\
p_1(r,\omega) &=& \Pr((r,\omega) \mbox{ is nearest city to $\mathbf{0}$ in $\cone(\mathbf{0},0,\pi/k)$}\\
&&\mbox{ and $\mathbf{0}$ is nearest city to $(r,\omega)$ in $\cone((r,\omega),\pi,\pi + \pi/k)$}) .
\end{eqnarray*}
We assert
\begin{equation}
L_k = \int_0^\infty \int_0^{\pi/k} r \ [2p(r,\omega) - p_1(r,\omega)] \  \ d\omega \ r \  dr . 
\label{Lpr}
\end{equation}
To argue this, first consider only roads $(v_L,v_R)$, written so that the $x$-coordinate of $v_L$ is less than the $x$-coordinate of $x_R$, and for which each city is the closest neighbor of the other city in the relevant cone.  
Given such cities at $v_L = \mathbf{0}, \ v_R = (r,\omega)$ the probability of such a road is $p_1(r,\omega)$ and the contribution to mean network
length is $r p_1(r,\omega)$.
Because the density of possible positions of $(v_L,v_R)$ has intensity $1$
on the region where $v_R \in \cone(v_L,0,\pi/k)$,
the contribution to
 normalized network length will be 
\[ \int_0^\infty \int_0^{\pi/k} r \ [ p_1(r,\omega) ]  \ \ d\omega \ r \  dr . \]
If instead we consider only roads $(v_L,v_R)$ where $v_R$ is the nearest neighbor to $v_L$ in its cone but not conversely, then similarly the normalized length of such roads is 
\[ \int_0^\infty \int_0^{\pi/k} r \ [ p (r,\omega) - p_1(r,\omega) ]  \ \ d\omega \ r \ dr . \]
By symmetry, the opposite possibility -- that $v_L$ is the nearest neighbor to $v_R$ in its cone but not conversely -- makes the same contribution.  Summing these three contributions gives (\ref{Lpr}).

To write formulas for $p(\cdot)$ and $p_1(\cdot)$, recall that the probability that the Poisson process assigns no cities to a region $A$ equals 
$\exp(- \area(A))$. 
For $p(\cdot)$, the relevant region is the finite cone $\mathbf{0} CE$ in Fig. 6, which has area 
$\sfrac{\pi r^2}{2k}$,
and so
\begin{equation}
p(r,\omega) = \exp(- \sfrac{\pi r^2}{2k}) .
\end{equation}
For $p_1(\cdot)$, the relevant region is the  entire region $\mathbf{0}ABCDEFG$ in Fig. 6.
The area of this region can be represented as 
\begin{quote}
area of cone $\mathbf{0} CE$, plus area of cone $DGA$, minus area of parallelogram $\mathbf{0}BDF$.
\end{quote}
The parallelogram  has height $r \sin \omega$ and base 
$r \cos \omega - r \frac{\sin \omega}{\tan \pi/k}$ and hence has area
\[ r^2  \ (\cos \omega -  \sfrac{\sin \omega}{\tan \pi/k}) \ \sin \omega. \]
So the area of $\mathbf{0}ABCDEFG$  equals 
\[ \sfrac{\pi r^2}{2k} + \sfrac{\pi r^2}{2k} - r^2  \ (\cos \omega -  \sfrac{\sin \omega}{\tan \pi/k}) \ \sin \omega \]
and finally 
\begin{equation}
p_1(r,\omega) = \exp( - r^2 [ \sfrac{\pi}{k} - \cos \omega \sin \omega + \sfrac{\sin^2 \omega}{\tan \pi/k}
] )
\end{equation}
Returning to formula (\ref{Lpr}), because
$\int_0^\infty r^2  \exp(-ar^2) \ dr = \sfrac{1}{4} \pi^{1/2} a^{-3/2}$, we can integrate out $r$ to get 
\begin{eqnarray}
L_k&=&
\sfrac{1}{4} \pi^{1/2} \int_0^{\pi/k}  
\left(2 (\sfrac{\pi}{2k})^{-3/2} \ - \ [ \sfrac{\pi}{k} - \cos \omega \sin \omega + \sfrac{\sin^2 \omega}{\tan \pi/k}
]^{-3/2} \right) \ d \omega \nonumber \\
&=& \sqrt{2k} - \sfrac{1}{4} \pi^{1/2} \int_0^{\pi/k}  
\ [ \sfrac{\pi}{k} - \cos \omega \sin \omega + \sfrac{\sin^2 \omega}{\tan \pi/k}
]^{-3/2}  \ d \omega 
\end{eqnarray}
which is formula (\ref{Lkformula}).

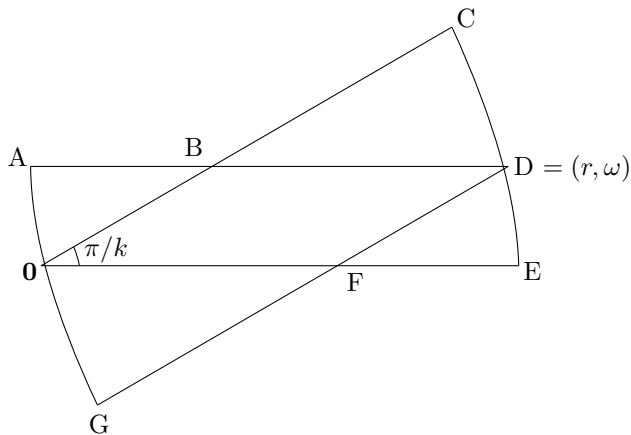
\begin{figure}
\setlength{\unitlength}{2.5in}
\begin{picture}(1.5,1.2)(-0.4,-0.4)
\put(0,0){\line(1,0){1}}
\put(0,0){\line(86,50){0.86}}
\put(0.978,0.208){\line(-1,0){1}}
\put(0.978,0.208){\line(-86,-50){0.86}}
\qbezier(1,0)(0.994,0.208)(0.86,0.5)
\qbezier(-0.022,0.208)(-0.02,0)(0.118,-0.292)
\put(-0.07,0.21){A}
\put(0.3,0.23){B}
\put(0.87,0.5){C}
\put(0.99,0.19){D $ = (r,\omega)$}
\put(1.01,-0.03){E}
\put(0.64,-0.052){F}
\put(0.1,-0.35){G}
\put(-0.04,-0.03){$\mathbf{0}$}
\put(0.09,0.02){$\pi/k$}
\qbezier(0.08,0)(0.077,0.019)(0.068,0.04)
\end{picture}
\caption{Regions of integration in the proof of Lemma \ref{LLk}.}
\end{figure}

\section{Lower bounds in long networks; average-case analysis}
\label{sec-6}
Turning to lower bounds, for $\Psave$ we start by giving  a reformulation 
(\ref{Poisson-interpret-2}) of the interpretation 
(\ref{Poisson-interpret}) in terms of a Poisson point process on the infinite plane. 
In (\ref{Poisson-interpret}) we required the distribution $\mu$ of the network to be translation invariant; by applying a 
random rotation $\Theta$ (uniform on $(0,2\pi)$) we may suppose also that $\mu$  is isotropic.  
Recall $L(\mu)$ and $S(\mu)$ denote normalized length and stretch.
Consider the number
\[ \intersect(\mu) = \mbox{ mean number of intersections of network edges with } \] \[ \mbox{ the $x$-axis per unit length.}
\]
There is a general formula 
(see \cite{MR895588} Chapter 8 for the relevant theory) that for any isotropic translation invariant network,
\begin{equation}
 L(\mu) = \sfrac{\pi}{2} \times \ \intersect(\mu).  \label{len-intersect}
\end{equation}
So we can rewrite (\ref{Poisson-interpret}) as
\begin{equation}
\Psave(s) = \sfrac{\pi}{2} \times \inf \{L(\mu); \ \mu \mbox{ is isotropic translation invariant, } S(\mu) \leq s \} .
\label{Poisson-interpret-2}
\end{equation}
We will use this formulation to obtain an order of magnitude lower bound for small $s$.  This general method was used in  somewhat different contexts in 
\cite{MR2411811,EJP2920}.
\begin{Proposition}
\label{P38}
$\Psave(s) = \Omega((s-1)^{-3/8})$ as $s \downarrow 1$.
\end{Proposition}
{\bf Proof.}
Given $h > 0$ consider the rate-$1$ Poisson point process restricted to the infinite strip 
$(-\infty, \infty) \times [-h,h]$.  
Consider pairs of such Poisson points, where one point is above the $x$-axis and the other is below the $x$-axis, and where the line segment between the two points crosses the $x$-axis at an angle greater that $45^\circ$. 
That is, consider pairs at positions 
$(x_1,y_1)$ and $(x_2,y_2)$ related by
\begin{equation}
-h < \min(y_1,y_2) < 0 < \max(y_1,y_2) < h; \quad 
|x_2-x_1| < |y_2 - y_1| 
 . \label{2-points}
\end{equation} 
Call such a pair {\em friends}.
For each friends pair, a hypothetical straight line segment between them crosses the $x$-axis at some position $\chi$, and the set of all such ``virtual crossing positions"  is a stationary point process on the line $(-\infty, \infty)$. 
For $L > 0$ write 
\[ N(h,L) = \mbox{number of virtual crossing positions in } [ 0,L] . \]

Now consider a network with stretch $\le 1+s$ over the rate-$1$ Poisson point process on the plane.  (So here $s>0$; this notational shift simplifies formulas).
The route between two friends
must cross the $x$-axis at some ``route-crossing position" $\chi^\prime$; write 
$\delta(h,s)$ for the maximum possible value of the distance between the 
route-crossing position $\chi^\prime$ and the virtual crossing position $\chi$.  
It is geometrically clear that this maximum is attained when the friends  are at positions 
$(-h,-h)$ and $(h,h)$, and therefore 
\begin{equation}
\delta(h,s) = h g^{-1}(s) 
\label{dhs}
\end{equation}
where $g^{-1}(\cdot)$ is the inverse function of 
\[
g(\delta) = \frac{\sqrt{1 + (1+\delta)^2} + \sqrt{1 + (1-\delta)^2} }{2 \sqrt{2}} - 1 \]
for which we calculate 
\begin{equation}
 g(\delta) \sim  \delta^2/8 \mbox{ as } \delta \downarrow 0 . 
\label{gd}
\end{equation}

Now choose $L > 0$ and partition the $x$-axis into blocks of length $L + 2 \delta(h,s)$, each block
consisting of a middle interval of length $L$ surrounded by two intervals of length $\delta(h,s)$. 
If the middle interval contains the  virtual crossing position for a pair of friends in the Poisson process,
then the block contains  the 
route-crossing position, 
and it follows that the rate of such  
route-crossing positions is at least
$\Pr(N(h,L) \ge 1) / (L + 2 \delta(h,s))$.  
We may choose $h$ and $L$ arbitrarily, so
appealing to (\ref{Poisson-interpret-2})  we have
\begin{equation}
 \Psave(1+s) \ge \sfrac{\pi}{2} \ \sup_{h,L} 
\frac{ \Pr(N(h,L) \ge 1) } { L + 2 \delta(h,s)} . 
\label{aveLB}
\end{equation}
We can lower bound the numerator via the second moment inequality 
\begin{equation}
 \Pr(N(h,L) \ge 1)  \ge 
\frac{ (\Ex N(h,L) )^2} { \Ex N^2(h,L)} . 
\label{2moment}
\end{equation}
It is eaasy to calculate $\Ex N(h,L) $, as follows.
 For a point $(0,-y_0)$ consider the set of possible positions 
of a friend $(x,y)$ with $x>0$.
The constraints are 
\[ 0 < y < h, \quad 0 < x < y_0 + y \]
and the area of this region equals 
$hy_0 + h^2/2$.
It follows easily that the rate of the stationary process of
virtual crossing positions equals 
\[ 2 \int_0^h (hy_0 + h^2/2) \ dy_0 = 2h^3 . \]
The initial factor $2$ arises due to the symmetric possibility 
$(0,+y_0)$ for the left point.
So we have shown 
\[ \Ex N(h,L) = 2h^3L   .\]
We will be concerned with the limit regime
\begin{equation}
h \to \infty, \ L \to 0, \ 2h^3L \to \lambda
\label{Nregime}
\end{equation}
for arbitrary $0 < \lambda < \infty$. 
Intuitively we  expect that the distribution of  $N(h,L)$ converges to Poisson$(\lambda)$ in this regime, but for our purposes it will suffice
to prove the second moment result (consistent with the Poisson limit)
\begin{equation}
\Ex N^2(h,L) \to \lambda^2 + \lambda 
\mbox{ in the limit regime } (\ref{Nregime}) .
\label{Nregime2}
\end{equation}
Defering the  proof of (\ref{Nregime2}), Proposition \ref{P38}
can be deduced from the ingredients above. Set
\[ h = h(s) = s^{-1/8}, \ L = L(s) = s^{3/8} \] 
and consider orders of magnitude as $s \downarrow 0$.  
The numerator in (\ref{aveLB}) is $\Omega(1)$ by (\ref{2moment}) and (\ref{Nregime2}).  And by (\ref{dhs})  and (\ref{gd}) we see that $\delta(h,s)$ 
is order $h s^{1/2} = s^{3/8}$, so the denominator in  (\ref{aveLB}) 
is order $s^{3/8}$, establishing the Proposition.

\medskip
\noindent
{\bf Proof of (\ref{Nregime2}).}
The formula for the second moment is given as (\ref{N2}) below.  
The term $\Ex N(h,L)$ arises from individual crossings, and the
 term $(\Ex N(h,L) )^2$ is the contribution from pairs of virtual crossing positions in $[0,L]$ for which the 4 end-points are all distinct.  
The integral term is the contribution from
 the case of two 
virtual crossing positions in $[0,L]$ with an end-point in common, say at 
$(x_0,-y_0)$ where $y_0>0$. 
This term involves
 the region 
$A(x_0,-y_0)$ containing the possible positions of a friend of $(x_0,-y_0)$
 for which the
virtual crossing position is in $[0,L]$. 
Fig. 7 shows this region, for a particular value of $(x_0,-y_0)$.
The integrand $\sfrac{1}{2}  (\mathrm{area}  \ A(x_0,-y_0))^2$ is the mean (conditioned on a point at $(x_0,-y_0)$) number of pairs of friends for which  both
virtual crossing positions (from friend to $(x_0,-y_0)$) are in $[0,L]$.
This  leads to the formula  
\begin{equation}
 \Ex N^2(h,L) = \Ex N(h,L)  +  (\Ex N(h,L) )^2 + 
2 \int \int_B  \sfrac{1}{2}  (\mathrm{area}  \ A(x_0,-y_0))^2 \ dx_0 dy_0 .
\label{N2}
\end{equation} 
We integrate over  the region $B$  of values for $(x_0,-y_0)$ which are consistent with a virtual crossing position in $[0,L]$.  This region $B$ can be decomposed as the union of four regions $B_0, B^\ell_1, B^r_1,B_2$ as shown in Fig. 8, 
wherein we are assuming $h> L/2$, which is true in the limit regime.

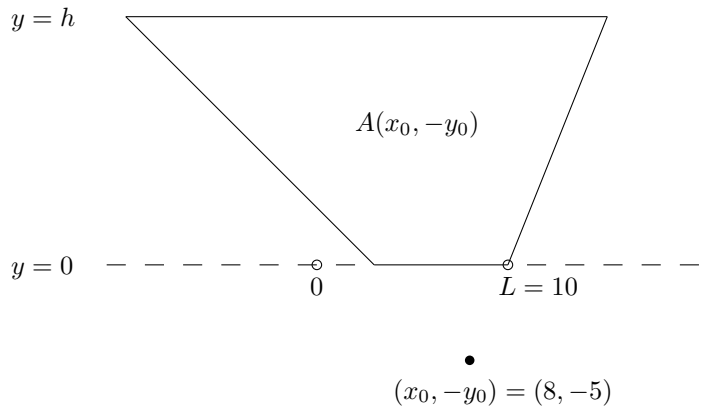
\begin{figure}
\setlength{\unitlength}{0.1in}
\begin{picture}(44,29)(-17,-9)
\put(3,0){\line(1,0){7}}
\put(3,0){\line(-1,1){13}}
\put(10,0){\line(2,5){5.2}}
\put(-10,13){\line(1,0){25.2}}
\put(-0.35,-1.6){0}
\put(9.5,-1.6){$L = 10$}
\put(8,-5){\circle*{0.5}}
\put(4,-7){$(x_0,-y_0) = (8,-5)$}
\multiput(11,0)(2,0){5}{\line(1,0){1}}
\multiput(-11,0)(2,0){7}{\line(1,0){1}}
\put(2,7){$A(x_0,-y_0)$}
\put(0,0){\circle{0.5}}
\put(10,0){\circle{0.5}}
\put(-16,-0.5){$ y = 0$}
\put(-16,12.5){$ y = h$}
\end{picture}
\caption{The region $A(x_0,-y_0)$ for the point $\bullet$.}
\end{figure}

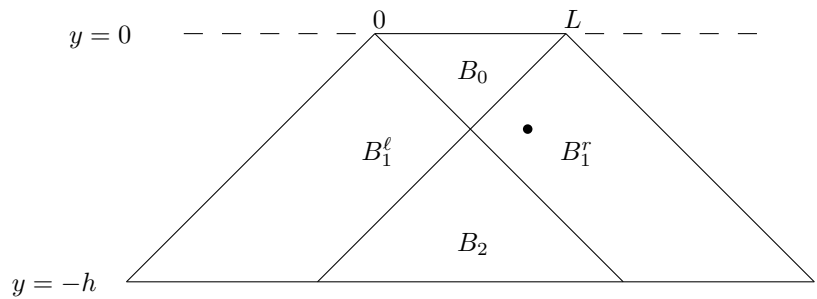
\begin{figure}
\setlength{\unitlength}{0.1in}
\begin{picture}(44,22)(-17,-17)
\put(0,0){\line(1,0){10}}
\put(0,0){\line(-1,-1){13}}
\put(10,0){\line(-1,-1){13}}
\put(0,0){\line(1,-1){13}}
\put(10,0){\line(1,-1){13}}
\put(-0.1,0.3){0}
\put(9.9,0.3){$L$}
\put(4.3,-2.4){$B_0$}
\put(9.7,-6.6){$B_1^r$}
\put(-0.7,-6.6){$B_1^\ell$}
\put(4.3,-11.4){$B_2$}
\put(-13,-13){\line(1,0){36}}
\multiput(11,0)(2,0){5}{\line(1,0){1}}
\multiput(-10,0)(2,0){5}{\line(1,0){1}}
\put(8,-5){\circle*{0.5}}
\put(-16,-0.5){$ y = 0$}
\put(-19,-13.5){$ y = - h$}
\end{picture}
\caption{The decomposition of the region of points consistent with a virtual crossing position in $[0,L]$.  The point $\bullet$ is the same as in Fig. 7.}
\end{figure}

For $(x_0,-y_0) \in B^r_1$, the case shown in Fig. 7, the 
region $A(x_0,-y_0)$ is the trapezoid bounded by the line 
$y = 0$, the line $y = h$, the line of slope $-1$ through $(x_0,-y_0)$ 
and the line through $(x_0,-y_0)$ and $(L,0)$. 
A brief calculation shows
\[ \mathrm{area}  \ A(x_0,-y_0) = \sfrac{1}{2} 
\left( 1 + \sfrac{L-x_0}{y_0} \right) \ 
\left( (h+y_0)^2 - y_0^2 \right) \quad \mbox{ for } 
(x_0,-y_0) \in B^r_1 .\] 
Easier calculations show 
\[ \mathrm{area}  \ A(x_0,-y_0) = 
 (h+y_0)^2 - y_0^2  \quad \mbox{ for } 
(x_0,-y_0) \in B_0 .\] 
\[ \mathrm{area}  \ A(x_0,-y_0) = \sfrac{L}{2 y_0} \ 
 \left((h+y_0)^2 - y_0^2 \right) \quad \mbox{ for } 
(x_0,-y_0) \in B_2 .\] 
The case $B^\ell_1$ is symmetric with $B^r_1$.
We could calculate  $\Ex N^2(h,L)$ exactly using (\ref{N2}), but we only need an upper bound.
The formulas above show that, as $x_0$ varies for fixed $y_0$, the quantity 
``$ \mathrm{area}  \ A(x_0,-y_0) $" takes its  maximum value on $B_0$ or $B_2$, and so
\[ \mathrm{area}  \ A(x_0,-y_0) \le  \left((h+y_0)^2 - y_0^2 \right) \ 
\min(1, \sfrac{L}{2 y_0}) . \]
So the integral term in (\ref{N2}) is bounded by
\[
\int_0^h (L + 2y_0) \ \left( 
 \left((h+y_0)^2 - y_0^2 \right) \ 
\min(1, \sfrac{L}{2 y_0}) 
\right)^2 \ dy_0 
. \]
The integral over $0 < y_0 < L/2$ works out as 
$\frac{3}{4} h^4L^2 + \frac{5}{6}h^3L^3 + \frac{7}{24}h^2L^4$.
The integral over $L/2 < y_0 < h$ works out as 
$\frac{7}{2} h^4L^2 - \frac{1}{4}h^3L^3 - \frac{3}{4}h^2L^4 
+ (\frac{1}{2}L^2h^4 + L^3h^3) \log (2h/L)$.  
So in the limit regime (\ref{Nregime}), the leading term is the term 
$\frac{1}{2}L^2h^4 \  \log (2h/L)$.  But this term $\to 0$, 
establishing (\ref{Nregime2}).

\section{Lower bounds on $\Psworst$ based on local optimality}
\label{sec-7}
One can get lower bounds on $\Psworst$ by choosing any configuration of cities and lower bounding the network length required for a network on that particular configuration to have a given stretch.
There are heuristic reasons (and the Steiner constant results mentioned at the start of
section \ref{sec-3}) to suspect that some kinds of regular configurations 
(rather than typical random configurations) are close to worst-case, so it is not unreasonable to use regular configurations to obtain lower bounds on worst-case behavior.
This allows us to work directly on the infinite plane, because the regular 
configurations we use have known average number of points per unit area.

\subsection{A bound from the square grid}
Consider, for instance, the ``square grid" configuration of cities at the points 
$\{(i,j); - \infty < i, j < \infty\}$.  
The usual ``square lattice" network (roads between city pairs $(v,w)$ at distance $1$)
has normalized length $= 2$ 
and stretch $= \sqrt{2}$. 
It is natural to conjecture this network is optimal amongst Steiner networks, in the following sense.
\begin{Conjecture}
If a Steiner network on the square grid configuration has stretch $\leq \sqrt{2}$ 
then its normalized length is at least $2$.
\end{Conjecture}
If true, this would imply 
$\Psworst(\sqrt{2} ) \geq 2$.
Similarly, any result of the type
\begin{quote}
A particular network $\NN_0$ on a particular configuration $\bz$ is optimal, in the sense that 
any other network $\NN$ with $\stretch(\NN) \leq \stretch(\NN_0) = s_0$ has normalized length 
$L(\NN) \geq L(\NN_0) = \ell_0$
\end{quote}
would imply $\Psworst(s_0) \geq \ell_0$.  
However, we are unable to prove {\em any} result of this type.  
Instead, we can only prove weaker results of the following type.
Consider the ``alternate diagonals" network on the square grid, shown in Fig. 9.

\begin{figure}
\setlength{\unitlength}{0.22in}
\begin{picture}(10,8)(-6,-1)
\multiput(0,0)(1,0){7}{\circle*{0.19}}
\multiput(0,0)(1,0){7}{\circle*{0.19}}
\multiput(0,1)(1,0){7}{\circle*{0.19}}
\multiput(0,2)(1,0){7}{\circle*{0.19}}
\multiput(0,3)(1,0){7}{\circle*{0.19}}
\multiput(0,4)(1,0){7}{\circle*{0.19}}
\multiput(0,5)(1,0){7}{\circle*{0.19}}
\multiput(0,6)(1,0){7}{\circle*{0.19}}
\put(1.5,-0.5){\line(1,1){5}}
\put(3.5,-0.5){\line(1,1){3}}
\put(5.5,-0.5){\line(1,1){1}}
\put(-0.5,-0.5){\line(1,1){7}}
\put(-0.5,1.5){\line(1,1){5}}
\put(-0.5,3.5){\line(1,1){3}}
\put(-0.5,5.5){\line(1,1){1}}
\put(-0.5,1.5){\line(1,-1){2}}
\put(-0.5,3.5){\line(1,-1){4}}
\put(-0.5,5.5){\line(1,-1){6}}
\put(0.5,6.5){\line(1,-1){6}}
\put(2.5,6.5){\line(1,-1){4}}
\put(4.5,6.5){\line(1,-1){2}}
\end{picture}
\caption{The ``alternate diagonals" network.}
\end{figure}
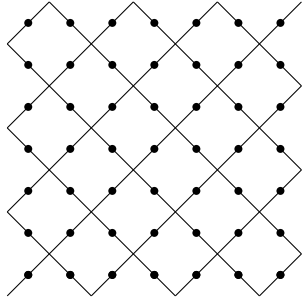

By inspection, this network has normalized length $= \sqrt{2}$ and satisfies
\begin{equation}
\mbox{route-length from $v$ to $w$ is $\leq \sqrt{2}$ for each city pair $(v,w)$ at Euclidean distance $1$}.
\label{diag}
\end{equation}
We can prove this network is optimal with respect to those properties.
\begin{Proposition}
\label{P11}
Any network on the square grid configuration satisfying (\ref{diag}) has normalized length 
$\geq \sqrt{2}$.
\end{Proposition}
\begin{Corollary}
\label{C11}
$\Psworst(\sqrt{2} ) \geq \sqrt{2}$.
\end{Corollary}
We call Proposition \ref{P11} a  ``local optimality" result because (\ref{diag}) is a ``local" analog of stretch.

\medskip
\noindent
{\bf Proof of Proposition \ref{P11}.}  
Take some network connecting the cities in the square grid configuration.
Consider a route through cities 
$\ldots \to (-2,0) \to (-1,0) \to (0,0) \to (1,0) \to (2,0) \to \ldots$ 
using minimum-length routes between each successive pair of cities. 
As we traverse this route, we might {\em backtrack}, meaning that the $x$-coordinate of position might decrease, but discarding any backtracking segments leaves a (maybe disconnected) 
non-backtracking route $((x,y(x)), - \infty < x < \infty)$.  Call this ``horizontal" route $H_0$.  Define the measure  $U_{H_0}$
on $H_0$ as the measured induced by Lebesgue measure on $x$; that is,
a line segment in $H_0$ from $(x_1,y(x_1))$ to $(x_2,y(x_2))$ has measure $x_2 - x_1$.
Repeat for routes $H_j$ through $\ldots \to (-2,j) \to (-1,j) \to (0,j) \to (1,j) \to (2,j) \to \ldots$.
The key observation is that assumption (\ref{diag}) implies that routes $H_j$ are disjoint as $j$ varies, except that routes $H_j$ and $H_{j+1}$ can meet at isolated points of the form $(i+\frac{1}{2},j+\frac{1}{2})$, as happens 
in the  ``alternate diagonals" network. 

 Let $\mu_H = \sum_{j= - \infty}^\infty \mu_{H_j}$.
It is clear that $\mu_H$ has ``density $1$", in the sense that for increasing  squares $A$
\begin{equation}
\frac{\mu_H(A)}{\area(A)} \to 1 \mbox{ as } \area (A) \to \infty . 
\label{muHV}
\end{equation}
Repeat the construction with vertical routes $V_i$ through $\ldots \to (i,-2) \to (i,-1) \to (i,0) \to (i,1) \to (i,2) \to \ldots$ to define a measure $\mu_V$ which also satisfies (\ref{muHV}).

Now write $\Lambda$ for length measure on the edges of the network. 
Consider a point $(x,y)$ on a road segment at angle $\theta$.
By the disjointness property, this point is in at most one $H_j$, in which case the density $d \mu_H/d\Lambda$ at the point equals 
$|\cos \theta|$, and in at most one $V_i$, 
in which case the density $d \mu_V/d\Lambda$ at the point equals 
$|\sin \theta|$.
It follows that 
\begin{equation}
 \frac{d(\mu_H + \mu_V)}{d \Lambda}(x,y) \leq |\cos \theta|  +  |\sin \theta| . 
\label{cos+sin}
\end{equation}
But always $|\cos \theta| + |\sin \theta| \leq \sqrt{2}$, so for any region $A$
\[ \mu_H(A) + \mu_V(A) \leq \sqrt{2}\  \Lambda(A) \] 
and then (\ref{muHV}) implies 
\[ \Lambda(A) \geq (\sqrt{2} - o(1)) \ \area(A)   \mbox{ as } \area (A) \to \infty . 
 \]
 That is, normalized network length is at least $\sqrt{2}$.

 \subsection{Another bound from hexagons}
 Here we show how the argument scheme above can be adapted to the hexagonal configuration of cities (Fig. 11).
 \begin{Proposition}
 \label{P12}
 Let $\NN$ be a network on the hexagonal configuration such that
 \begin{equation}
 \frac{
r(v,w)}
{d(v,w)} \leq \sqrt{3} 
\mbox{ for all (Euclidean) nearest-neighbor pairs } (v,w) . 
\label{hex-ass}
\end{equation}
Then its normalized length is at least 
$2^{-1} 3^{3/4}$. 
\end{Proposition}
 \begin{Corollary}
 \label{C12}
$\Psworst(\sqrt{3} ) \geq 2^{-1} 3^{3/4} =  1.14 ....$.
\end{Corollary}
 {\bf Proof of Proposition \ref{P12}.}  
 Consider the hexagonal configuration with $\ell = $ distance between nearest neighbors.
 The density of cities (number per unit area) is 
 \begin{equation}
  \rho(\ell) = 4 \cdot 3^{-3/2} \ \ell^{-2} . 
 \label{rho-ell}
\end{equation}
 Fig. 10 shows four adjacent cities $ABCD$ in one hexagon.  
 In that figure we see the route lengths satisfy
 \[ \frac{\len (AZB)}{d(A,B)} = \frac{\len(DZC)}{d(D,C)}  = \sqrt{3} \]
 and it is easy to check the optimality property:
 \begin{quote}
 if $\pi_1$ and $\pi_2$ are paths in the plane from $A$ to $B$ and from $C$ to $D$ respectively, 
 and if 
 $\max( \frac{\len (\pi_1)}{d(A,B)} , \frac{\len(\pi_2)}{d(D,C)} )\leq \sqrt{3}$, 
 then the paths cannot meet except possibly at $Z$.
 \end{quote}

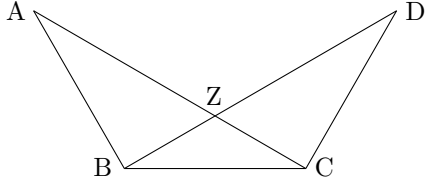
\begin{figure}
 \setlength{\unitlength}{0.95in}
\begin{picture}(1.5,1.4)(-0.7,-0.2)
\put(-0.5,0){\line(1,0){1}}
\put(0.5,0){\line(50,87){0.5}}
\put(-0.5,0){\line(-50,87){0.5}}
\put(-0.5,0){\line(150,87){1.5}}
\put(0.5,0){\line(-150,87){1.5}}
\put(-1.15,0.82){A}
\put(-0.67,-0.04){B}
\put(0.55,-0.04){C}
\put(1.05,0.82){D}
\put(-0.05,0.35){Z}
 \end{picture}
\caption{An optimality property.}
 \end{figure}

 Now consider the minimum-length route in $\NN$ through an ``angle $\pi/6$ staircase" like 
 $abcdef \ldots$ in Fig. 11.  
 By assumption (\ref{hex-ass}) and the optimality property above, this route does not meet the corresponding route through the next staircase $ghijkl \ldots$ except at isolated points.
 As in the previous section, each path segment on such a route is at some angle $\theta$ to the 
 ``angle $=\pi/6$" line; put a measure on the non-backtracking parts of the route with density 
 $\cos \theta$ w.r.t. length measure $\Lambda$ on the segment.  
 Repeating for each  angle $=\pi/6$ staircase gives a measure $\mu_{\pi/6}$ on network edges, 
 which has the property (for squares $A$)
 \[
\frac{\mu_{\pi/6}(A)}{\lambda_{\pi/6}(A)} \to 1 \mbox{ as } \area (A) \to \infty 
\]
where $\lambda_{\pi/6}$ is length measure on the parallel  ``angle $=\pi/6$" {\em straight} lines through the staircases.
The orthogonal distance between such lines equals $3 \ell/2$ 
(this is easiest to see with the angle $=\pi/2$ lines, where the distance is the average of  $d(c,d)$ and $d(b,k)$), 
so
\[
\frac{\lambda_{\pi/6}(A)}{\area(A)} \to \frac{2}{3 \ell}  \mbox{ as } \area (A) \to \infty 
\]
and thus
\begin{equation}
\frac{\mu_{\pi/6}(A)}{\area(A)} \to \frac{2}{3 \ell}  \mbox{ as } \area (A) \to \infty . 
\label{pi6}
\end{equation}

 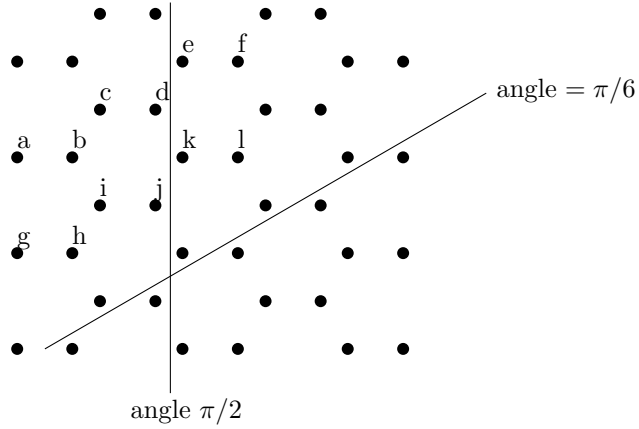
\begin{figure}
 \setlength{\unitlength}{0.33in}
 \begin{picture}(10,8)(7,18)
 \multiput(10,20)(0.8775,0){2}{\circle*{0.19}}
\multiput(12.6225,20)(0.8775,0){2}{\circle*{0.19}}
\multiput(15.245,20)(0.8775,0){2}{\circle*{0.19}}
\multiput(11.316,20.76)(0.8775,0){2}{\circle*{0.19}}
\multiput(13.9385,20.76)(0.8775,0){2}{\circle*{0.19}}
\multiput(10,21.52)(0.8775,0){2}{\circle*{0.19}}
\multiput(12.6225,21.52)(0.8775,0){2}{\circle*{0.19}}
\multiput(15.245,21.52)(0.8775,0){2}{\circle*{0.19}}
\multiput(11.316,22.28)(0.8775,0){2}{\circle*{0.19}}
\multiput(13.9385,22.28)(0.8775,0){2}{\circle*{0.19}}
\multiput(10,23.04)(0.8775,0){2}{\circle*{0.19}}
\put(10.0,23.19){a}
\put(10.8775,23.19){b}
\put(11.316,23.95){c}
\put(12.193,23.95){d}
\put(12.6225,24.71){e}
\put(13.5,24.71){f}
\put(10.0,21.67){g}
\put(10.8775,21.67){h}
\put(11.316,22.43){i}
\put(12.193,22.43){j}
\put(12.6225,23.19){k}
\put(13.5,23.19){l}
\multiput(12.6225,23.04)(0.8775,0){2}{\circle*{0.19}}
\multiput(15.245,23.04)(0.8775,0){2}{\circle*{0.19}}
\multiput(11.316,23.8)(0.8775,0){2}{\circle*{0.19}}
\multiput(13.9385,23.8)(0.8775,0){2}{\circle*{0.19}}
\multiput(10,24.56)(0.8775,0){2}{\circle*{0.19}}
\multiput(12.6225,24.56)(0.8775,0){2}{\circle*{0.19}}
\multiput(15.245,24.56)(0.8775,0){2}{\circle*{0.19}}
\multiput(11.316,25.32)(0.8775,0){2}{\circle*{0.19}}
\multiput(13.9385,25.32)(0.8775,0){2}{\circle*{0.19}}
\put(10.4388,20){\line(150,87){7}}
\put(17.6,24){angle $= \pi/6$}
\put(12.43,19.3){\line(0,1){6.2}}
\put(11.8,18.9){angle $\pi/2$}
\end{picture}
\caption{The hexagonal configuration.  Points $abcdef$ are on a 
``angle $= \pi/6$ staircase" parallel (in an symptotic sense) to the 
``angle $= \pi/6$ line".
}
\end{figure}

\vspace{0.1in}
\noindent
Repeat the construction with staircases like $cdkl \ldots$ with angle $= -\pi/6$ to get a measure $\mu_{-\pi/6}$ on the associated routes; repeat again with staircases like $jkde\ldots$ with angle $= -\pi/2$ to get a measure $\mu_{\pi/2}$.  
These measures also satisfy (\ref{pi6}). 
Note each adjacent pair of cities is in two staircases, of different angles. 
The analog of (\ref{cos+sin}) is that, 
at a point $(x,y)$ on a road segment at angle $\theta$,
\[
 \frac{d(\mu_{\pi/2} + \mu_{\pi/6} + \mu_{-\pi/6})}{d \Lambda}(x,y) \leq |\cos (\theta - \pi/2)|  +  |\cos (\theta - \pi/6)|  + |\cos (\theta + \pi/6)| 
\]
because the point is in at most one route for each of the three angles. 
But 
\[ |\cos (\theta - \pi/2)|  +  |\cos (\theta - \pi/6)|  + |\cos (\theta + \pi/6)|  \leq 2 \]
and so 
\[ (\mu_{\pi/2} + \mu_{\pi/6} + \mu_{-\pi/6})(A) \leq 2 \Lambda(A). \]
Use (\ref{pi6}) to see 
\begin{equation}
 \Lambda(A) \geq (\ell^{-1} - o(1)) \ \area(A)   \mbox{ as } \area (A) \to \infty . 
 \label{LAA}
 \end{equation}
 Our normalization convention is that cities have density $1$, that is $\rho(\ell) = 1$ at 
 (\ref{rho-ell}), so $\ell = 2 \cdot 3^{-3/4}$ and the lower bound in 
(\ref{LAA}) becomes 
 $\ell^{-1} = 2^{-1} 3^{3/4}$.

\subsection{The triangular lattice}
We sketch the minor modification which uses the triangular lattice (Fig. 12).

 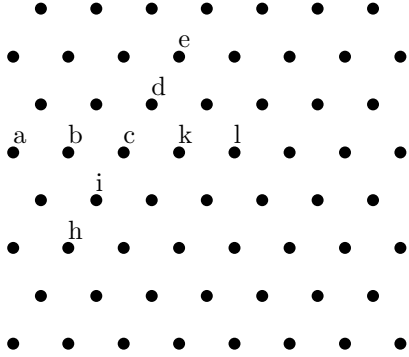
\begin{figure}
 \setlength{\unitlength}{0.33in}
 \begin{picture}(10,6.8)(7,19.5)
 \multiput(10,20)(0.8775,0){8}{\circle*{0.19}}
\multiput(10.441,20.76)(0.8775,0){7}{\circle*{0.19}}
\multiput(10,21.52)(0.8775,0){8}{\circle*{0.19}}
\multiput(10.441,22.28)(0.8775,0){7}{\circle*{0.19}}
\multiput(10,23.04)(0.8775,0){8}{\circle*{0.19}}
\put(10.0,23.19){a}
\put(10.8775,23.19){b}
\put(11.75,23.19){c}
\put(12.193,23.95){d}
\put(12.6225,24.71){e}
\put(10.8775,21.67){h}
\put(11.316,22.43){i}
\put(12.6225,23.19){k}
\put(13.5,23.19){l}
\multiput(10.441,23.8)(0.8775,0){7}{\circle*{0.19}}
\multiput(10,24.56)(0.8775,0){8}{\circle*{0.19}}
\multiput(10.441,25.32)(0.8775,0){7}{\circle*{0.19}}
\end{picture}
\caption{The triangular configuration.}
\end{figure}

 \begin{Proposition}
 \label{P13}
 Let $\NN$ be a network on the triangular configuration such that
 \begin{equation}
 \frac{r(v,w)}
{d(v,w)} \leq  
\sfrac{1}{2} + \sqrt{\sfrac{3}{4}} 
\mbox{ for all (Euclidean) nearest-neighbor pairs } (v,w) . 
\label{hex-ass}
\end{equation}
Then its normalized length is at least 
$2^{-1/2} 3^{3/4}$. 
\end{Proposition}
 \begin{Corollary}
 \label{C13}
$\Psworst \left(\sqrt{\sfrac{3}{4}} + \sfrac{1}{2} \right) \geq 2^{-1/2} 3^{3/4} =  1.61 ....$.
\end{Corollary}
 {\bf Outline proof of Proposition \ref{P13}.}  
We indicate changes in the previous argument.
The density of cities is now 
 \begin{equation}
  \rho(\ell) = 2 \cdot 3^{-1/2} \ \ell^{-2} . 
 \label{rho-ell-2}
\end{equation}
 Fig. 13 shows four adjacent cities $ABCD$ in  
the triangular configuration.
 In that figure we see
 \[ \frac{\len (AZD)}{d(A,D)} = \frac{\len(BZD)}{d(B.C)}   = 
\frac{1}{2} + \sqrt{\frac{3}{4}}
\]
 and it is easy to check the optimality property:
 \begin{quote}
 if $\pi_1$ and $\pi_2$ are paths in the plane from $A$ to $D$ and from $C$ to $B$ respectively, 
 and if 
 $\max( \frac{\len (\pi_1)}{d(A,D)} , \frac{\len(\pi_2)}{d(B,C)} )\leq  
\frac{1}{2} + \sqrt{\frac{3}{4}}
$,
 then the paths cannot meet except possibly at $Z$.
 \end{quote}

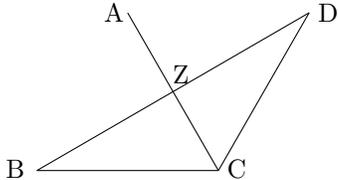
\begin{figure}
 \setlength{\unitlength}{0.95in}
\begin{picture}(1.5,1.4)(-0.4,-0.2)
\put(-0.5,0){\line(1,0){1}}
\put(0.5,0){\line(50,87){0.5}}
\put(-0.5,0){\line(150,87){1.5}}
\put(0.5,0){\line(-50,87){0.5}}
\put(-0.13,0.82){A}
\put(-0.67,-0.04){B}
\put(0.55,-0.04){C}
\put(1.05,0.82){D}
\put(0.25,0.48){Z}
 \end{picture}
\caption{An optimality property.}
 \end{figure}

As before,
there is a measure $\mu_0$ on routes through cities like 
$abckl$
on ``angle $ = 0$" routes,
and measures $\mu_{\pi/3}$ and $\mu_{-\pi/3}$ 
associated with angles $\pi_3$ (like $hicde$) and $-\pi/3$.
These satisfy 
\[
 \frac{d(\mu_{0} + \mu_{\pi/3} + \mu_{-\pi/3})}{d \Lambda}(x,y) \leq |\cos (\theta )|  +  |\cos (\theta - \pi/3)|  + |\cos (\theta + \pi/3)| 
\leq 2 .
\]
The orthogonal distance between parallel lines is 
$\ell \sqrt{3/4}$,
and repeating the argument for (\ref{LAA})  leads to 
\[ \Lambda(A) \geq (3^{1/2}\ell^{-1} - o(1)) \ \area(A)   \mbox{ as } \area (A) \to \infty . 
 \]
Taking $\rho(\ell) = 1$ in (\ref{rho-ell-2}),  the lower bound on normalized length is 
$3^{1/2} \ell^{-1} = 
2^{-1/2} 3^{3/4}$.

\subsection{Other configurations}
\label{sec-vineyard}
One could seek to repeat the arguments above with less symmetric configurations, but the calculations become messier, and we have not pursued details.

\paragraph{Acknowledgments}
We thank Lisha Li for  careful reading of a penultimate draft. 
Research supported  by N.S.F Grant DMS-1106998.


\begin{thebibliography}{10}

\bibitem{EJP2920}
David Aldous.
\newblock Scale-invariant random spatial networks.
\newblock {\em Electron. J. Probab.}, 19:no. 15, 1--41, 2014.

\bibitem{MR2411811}
David~J. Aldous and Wilfrid~S. Kendall.
\newblock Short-length routes in low-cost networks via {P}oisson line patterns.
\newblock {\em Adv. in Appl. Probab.}, 40(1):1--21, 2008.

\bibitem{MR2235175}
David~J. Aldous and Maxim Krikun.
\newblock Percolating paths through random points.
\newblock {\em ALEA Lat. Am. J. Probab. Math. Stat.}, 1:89--109, 2006.

\bibitem{MR2791668}
David~J. Aldous and Julian Shun.
\newblock Connected spatial networks over random points and a route-length
  statistic.
\newblock {\em Statist. Sci.}, 25(3):275--288, 2010.

\bibitem{bose_optimal}
Prosenjit Bose, Jean-Lou~De Carufel, Pat Morin, Andr{\'e} van Renssen, and
  Sander Verdonschot.
\newblock Optimal bounds on theta-graphs: More is not always better.
\newblock In {\em Proceedings of the 24th Canadian Conference on Computational
  Geometry, CCCG 2012, Charlottetown, Prince Edward Island, Canada, August
  8-10, 2012}, pages 291--296, 2012.

\bibitem{bose2013}
Prosenjit Bose, Andr\'{e} Renssen, and Sander Verdonschot.
\newblock On the spanning ratio of theta-graphs.
\newblock In {\em Algorithms and Data Structures}, Lecture Notes in Computer
  Science, pages 182--194. Springer, 2013.

\bibitem{MR627537}
F.~R.~K. Chung and R.~L. Graham.
\newblock On {S}teiner trees for bounded point sets.
\newblock {\em Geom. Dedicata}, 11(3):353--361, 1981.

\bibitem{abs-1305-4170}
Vida Dujmovic, Pat Morin, and Michiel H.~M. Smid.
\newblock Average stretch factor: How low does it go?, 2013.
\newblock http://arxiv.org/abs/1305.4170.

\bibitem{MR1134449}
J.~Mark Keil and Carl~A. Gutwin.
\newblock Classes of graphs which approximate the complete {E}uclidean graph.
\newblock {\em Discrete Comput. Geom.}, 7(1):13--28, 1992.

\bibitem{lando}
Tamar Lando.
\newblock Efficient networks and enumerations on forests.
\newblock Master's thesis, U.C. Berkeley, 2008.
\newblock http://www.stat.berkeley.edu/$\sim$aldous/Papers/Lando\_thesis.pdf.

\bibitem{MR0279853}
R.~E. Miles.
\newblock On the homogeneous planar {P}oisson point process.
\newblock {\em Math. Biosci.}, 6:85--127, 1970.

\bibitem{morin_v}
Pat Morin and Sander Verdonschot.
\newblock On the average number of edges in theta graphs, 2013.
\newblock http://arxiv.org/abs/1304.3402.

\bibitem{MR2289615}
Giri Narasimhan and Michiel Smid.
\newblock {\em Geometric spanner networks}.
\newblock Cambridge University Press, Cambridge, 2007.

\bibitem{MR805539}
Franco~P. Preparata and Michael~Ian Shamos.
\newblock {\em Computational geometry}.
\newblock Texts and Monographs in Computer Science. Springer-Verlag, New York,
  1985.
\newblock An introduction.

\bibitem{MR1422018}
J.~Michael Steele.
\newblock {\em Probability theory and combinatorial optimization}, volume~69 of
  {\em CBMS-NSF Regional Conference Series in Applied Mathematics}.
\newblock Society for Industrial and Applied Mathematics (SIAM), Philadelphia,
  PA, 1997.

\bibitem{MR895588}
D.~Stoyan, W.~S. Kendall, and J.~Mecke.
\newblock {\em Stochastic geometry and its applications}.
\newblock Wiley Series in Probability and Mathematical Statistics: Applied
  Probability and Statistics. John Wiley \& Sons Ltd., Chichester, 1987.
\newblock With a foreword by D. G. Kendall.

\bibitem{yao1982}
A.~Yao.
\newblock On constructing minimum spanning trees in k-dimensional spaces and
  related problems.
\newblock {\em SIAM Journal on Computing}, 11(4):721--736, 1982.

\bibitem{MR1632875}
Joseph~E. Yukich.
\newblock {\em Probability theory of classical {E}uclidean optimization
  problems}, volume 1675 of {\em Lecture Notes in Mathematics}.
\newblock Springer-Verlag, Berlin, 1998.

\end{thebibliography}

\end{document}